\documentclass[DIV=14,letterpaper]{scrartcl}

\usepackage{lineno,hyperref}
\modulolinenumbers[5]

\usepackage{amsmath,amssymb,amsfonts,amsthm,subfig} 

\usepackage{graphicx}
\def \D{\mathbb{D}}
\def \DK{\mathbb{D}_K}
\def \M{\mathbb{M}}
\def \MK{\mathbb{M}_K}
\def \Th{\mathcal{T}_h}
\newcommand{\V}[1]{\mbox{\boldmath $ #1 $}}

\newcommand{\bey}{\begin{eqnarray}}
\newcommand{\eey}{\end{eqnarray}}
\newcommand{\nn}{\nonumber}
\newcommand{\beq}{\begin{equation}}
\newcommand{\eeq}{\end{equation}}
\theoremstyle{plain}

\theoremstyle{definition}

\theoremstyle{remark}
\newtheorem{exam}{\hspace{6mm}{\textbf Example}}[section]
\newtheorem{rem}{\hspace{6mm}{\textbf Remark}}[section]

\begin{document}

\date{}
		
\title{Anisotropic Mesh Adaptation for Finite Element Solution of Anisotropic Porous Medium Equation}
\author{Xianping Li%
\thanks{Department of Mathematics and Statistics, the University of Missouri-Kansas City, Kansas City, MO 64110, U.S.A. (\textit{lixianp@umkc.edu}) }
}

\maketitle

\vspace{10pt}

\begin{abstract}
Anisotropic Porous Medium Equation (APME) is developed as an extension of the Porous Medium Equation (PME) for anisotropic porous media. A special analytical solution is derived for APME for time-independent diffusion. Anisotropic mesh adaptation for linear finite element solution of APME is discussed and numerical results for two dimensional examples are presented. The solution errors using anisotropic adaptive meshes show second order convergence. 
\end{abstract}

\noindent
\textbf{AMS 2010 Mathematics Subject Classification.} 65M60, 65M50

\noindent
\textbf{Key words.} {porous medium equation, anisotropic mesh adaptation, adaptive mesh, anisotropic diffusion, finite element, moving mesh}

\section{Introduction}
\label{Sec-intro}

In this paper, we extend the porous medium equation (PME) to the Anisotropic Porous Medium Equation (APME) that takes into consideration the anisotropic physical properties of the porous media such as permeability. Then we study the linear finite element solution of APME. We consider the following problem
\begin{equation}
\label{ibvp-1}
\begin{cases}
u_t = \nabla \cdot (u^m \D \nabla u), & \quad \mbox{ in } \quad \Omega_T = \Omega \times (t_0,T]
\\
u(\V{x},t) = 0, & \quad \mbox{ on } \quad \partial \Omega \times (t_0, T], \\
u(\V{x},0) = u_0(\V{x}), & \quad \mbox{ in } \quad \Omega \times \{t=t_0\}
\end{cases}
\end{equation}
where $u=u(\V{x},t)$ is a nonnegative scalar function, 
$m \ge 1$ is the physical parameter, $\Omega \subset \mathbb{R}^d$ is a connected polygonal or polyhedral domain of $d$-dimensional space,
$t_0 \ge 0$ is the starting time, $T > 0$ is the end time, $u_0(\V{x}) \ge 0$ is a given function.
We assume that $\mathbb{D}= \mathbb{D}(\V{x}) $ is a general symmetric and strictly positive
definite matrix-valued function on $\Omega_T$ that takes both isotropic and anisotropic diffusion
as special cases. For simplicity, we consider only time-independent diffusion matrix $\mathbb{D}$ in this work. 
The principles can also be applied to the time-dependent situation with minor modifications.

Porous medium equation (PME) arises in many fields of science and engineering such as fluid flow in porous media, heat transfer or diffusion, image processing, and population dynamics \cite{Vaz07}. The general form is given as
\begin{equation}
\label{pme-1}
u_t = \nabla \cdot (\nabla u^{m+1}),
\end{equation}
or in the modified form
\begin{equation}
\label{pme-2}
u_t = \nabla \cdot (u^{m}\nabla u).
\end{equation}
For gas flow in porous media, $m$ is the heat capacity ratio, $u$ represents the density, $u^{m}$ is the pressure, and $-\nabla u^{m}$ is the velocity.

Mathematically, the parameter $m$ in \eqref{pme-2} can take any real value. Specifically, when $m=0$, the PME \eqref{pme-1} or \eqref{pme-2} reduces to the heat equation. When $m>0$, the PME becomes a nonlinear evolution equation of parabolic type that has attracted interests of both theoretical and computational mathematicians. Particularly, when $m=1$, the PME is called the Boussinesq's equation that models groundwater flow in a porous stratum.

The nonlinear term $u^{m}$ in \eqref{pme-2} induces the so-called nonlinear diffusion that brings up many challenges in the analysis of the PME. However, this nonlinear diffusion is not the physical property of the porous media such as permeability. In this paper, we generalize PME to Anisotropic Porous Medium Equation (APME) that takes into account the anisotropy of the physical property of the porous media. In the mean time, APME can also be viewed as a nonlinear anisotropic diffusion problem. 

General PME \eqref{pme-1} or \eqref{pme-2} has been studied extensively both in theory \cite{OKC58, Aro69, Kal74, GP76, GP77, ACK83, Shm03, Vaz07} and numerical approximations \cite{Ros83, DH84, BCHR99, WL99, EL08, ZW09, PSRV10, ES12, SB12, DAA13, DAA15, LBL15, NH17}. In particular, the nonlinear diffusion term and the sharp gradient near the free boundary make it difficult to achieve high order convergence of the numerical solutions. For example, finite difference moving mesh method has been developed in \cite{SB12, LBL15} for PME in one-dimensional space and second order convergence has been observed; however, no result is provided for PME in 2D. Error estimates have been developed in \cite{PSRV10} for finite volume discretization of PME in 2D, which shows only first order convergence. For finite element discretization using quasi-uniform meshes \cite{Ros83, NV88, RW96, Ebm98, WL99}, the convergence is at most first order for $m=1$ and decreases for larger $m$. For one dimensional PME, high order convergence rate was achieved on a uniform mesh by using a high order local discontinuous Galerkin finite element method \cite{ZW09}. 

On the other hand, adaptive meshes and moving meshes are great choices to improve computational efficiency and accuracy by concentrating more mesh elements in the regions where solution changes significantly. In particular, Ngo and Huang \cite{NH17} have studied moving mesh finite element solution of PME and demonstrated the advantages of moving mesh over quasi-uniform meshes. Their results show second-order convergence of solution errors using moving meshes.  

However, no result is currently available for APME. Interesting features of PME also appear in APME such as finite propagation, free boundaries and waiting time phenomenon. Moreover, with the anisotropy of the porous media, satisfaction of maximum principle becomes more challenging and special mesh adaptation is needed, see \cite{KSS09, LH10, LH13} and the references therein. 
This paper serves as a starting effort about APME and its numerical solutions. Anisotropic mesh adaptation technique is applied in the numerical computations to improve efficiency and accuracy. Different than moving mesh method that keeps the connectivity of the elements, our anisotropic mesh adaptation technique can change the connectivity as well as the number of elements as desired. Therefore, a coarse uniform mesh can be used as the initial mesh in our adaptation, while a fine initial mesh is usually needed for moving mesh method in order to capture the sharp change of solution on the initial free boundary.

The outline of this paper is as follows. In Section \ref{Sec-APME}, the anisotropic porous medium equation (APME) is developed and the exact solution for a special case is discussed. Section \ref{Sec-fem} gives a brief summary of linear finite element solution of the APME, and Section \ref{Sec-adap} introduces the anisotropic mesh adaptation methods. Numerical examples are presented in Section \ref{Sec-results} to show the different behavior between APME and PME and also demonstrate the advantages of adaptive anisotropic meshes over other meshes. Finally, some conclusions are drawn in Section \ref{Sec-summary}.

\section{Anisotropic Porous Medium Equation (APME)}
\label{Sec-APME}
Firstly, we derive the model for fluid flow through anisotropic porous media as follows. Let $u$ be the density of the fluid. The flow is governed by the following three equations. 

(I) Conservation of mass (continuity equation)
\begin{equation}
\label{eq-mass}
\varepsilon u_t = - \nabla \cdot (u \, \V{v}),
\end{equation}
where $\varepsilon \in (0,1)$ is the porosity of the media and $\V{v}$ is the velocity. 

(II) Darcy's law in anisotropic porous media
\begin{equation}
\label{eq-velocity}
\V{v} = - \frac{1}{\mu} \D \nabla p,
\end{equation}
where $\mu$ is the viscosity of the fluid, $\D$ is the permeability matrix of the porous media, and $p$ is the pressure. In most cases, the porous media is anisotropic, thus $\D$ has different eigenvalues. If $\D$ varies with location, then it also represents heterogeneous media.

(III) The equation of state
\begin{equation}
\label{eq-state}
p= p_0 u^{m},
\end{equation}
where $p_0$ is the reference pressure and $m \ge 1$ is the ratio of specific heats. 

Combine \eqref{eq-mass}, \eqref{eq-velocity}, \eqref{eq-state} together, we have
\bey
\nn
u_t & = & \frac{p_0}{\varepsilon \mu} \nabla \cdot (u \, \D \nabla u^{m}) \\
\nn
& = & \frac{p_0}{\varepsilon \mu} \nabla \cdot (\D \, u \cdot m \cdot u^{m-1} \nabla u) \\
\label{apme0}
& = & \frac{m p_0}{\varepsilon \mu} \nabla \cdot (\D \, u^{m} \nabla u) \\
\label{apme1}
& = & \frac{m p_0}{(m+1) \varepsilon \mu} \nabla \cdot (\D \nabla u^{m+1}).
\eey
Scale time $t$ with the constant $\frac{m p_0}{(m+1) \varepsilon \mu}$ in \eqref{apme1} or $\frac{m p_0}{\varepsilon \mu}$ in \eqref{apme0}, we obtain the anisotropic porous medium equation (APME) as follows
\begin{equation}
\label{apme2}
u_t = \nabla \cdot (\D \nabla u^{m+1})	
\end{equation}
or
\begin{equation}
\label{apme}
u_t = \nabla \cdot (u^{m} \D \nabla u).	
\end{equation}
Combining $u^{m} \D$ in \eqref{apme} together as the diffusion term, APME can be viewed as an anisotropic diffusion equation with induced nonlinearity from the solution. Thus, available results for anisotropic diffusion problems \cite{BKK08, TW08, KSS09, LH10, LH13} are useful for the study of APME. 

Note that, V{\'a}zquez has mentioned the extension of PME in non-homogeneous media (NHPME) in \cite{Vaz07} as
\begin{equation}
\label{nhpme}
\varepsilon(\V{x},t)u_t = \nabla \cdot (c(\V{x},t)\nabla u^{m+1}),
\end{equation}
where $\varepsilon(\V{x},t)$ and $c(\V{x},t)$ are nonnegative functions. However, the anisotropy of the porous media has not been considered in NHPME. In this sense, both PME \eqref{pme-1} and NHPME \eqref{nhpme} are special cases of APME \eqref{apme}, where both the anisotropy and heterogeneity are taken into account in the diffusion matrix $\D$.  

Next, we derive a special solution for \eqref{ibvp-1}. 
We start with the Barenblatt-Pattle solution for PME \eqref{pme-1} developed by Barenblatt \cite{Bar52} and Pattle \cite{Pat59} independently. 
Suppose the initial solution $u_0(\V{x})$ (at time $t_0$) is compact-supported in a region given by $r_0 > 0$ as
\begin{equation}
\label{soln-pme-ini}
u_0(\V{x}) = \max \left\{ \left( 1 - \frac{\V{x}^T \V{x}}{r_0^2} \right)^{\frac{1}{m}}, \quad 0 \right\}.
\end{equation}
The solution at time $t \ge t_0$ for PME is given by 
\begin{equation}
\label{soln-pme}
u(\V{x},t) = \max \left\{ \frac{1}{\kappa^d} \left( 1 - \frac{\V{x}^T \V{x}}{r_0^2 \kappa^2} \right)^{\frac{1}{m}}, \quad 0 \right\},
\end{equation}
where 
\begin{equation}
\label{soln-para}
\kappa=\left(\frac{t}{t_0} \right)^{\beta}, \quad \beta=\frac{1}{dm+2}, \quad \text{ and } \quad  t_0=\frac{1}{2} \beta m r_0^2. 
\end{equation}

For APME, we consider the mesh $\Th$ in the physical domain $\Omega$ as a uniform mesh $\mathcal{T}_c$ in the computational domain $\Omega_c$ specified by the metric $\M$. By choosing $\M=\D^{-1}$, the diffusion can be considered as isotropic in the computational domain $\Omega_c$ with mesh $\mathcal{T}_c$ \cite{Hua06, LH10}. We apply the Barenblatt-Pattle solution for PME in $\Omega_c$ and then map the solution to $\Omega$ for APME. Denote the vertices in $\mathcal{T}_c$ by $\tilde{\V{x}}$, we have
\begin{equation}
	\tilde{\V{x}} = \M^{\frac{1}{2}} \V{x} = \D^{-\frac{1}{2}} \V{x}.
\end{equation}
We also consider the initial solution in a region given by $r_0 > 0$ in $\Omega_c$ as
\begin{equation}
\label{soln-apme-ini}
u_0(\V{x}) = \max \left\{ \left( 1 - \frac{\tilde{\V{x}}^T \tilde{\V{x}}}{r_0^2} \right)^{\frac{1}{m}}, \quad 0 \right\}
= \max \left\{ \left( 1 - \frac{\V{x}^T \D^{-1} \V{x}}{r_0^2} \right)^{\frac{1}{m}}, \quad 0 \right\}.
\end{equation}
Then the solution at time $t \ge t_0$ for APME is given by
\begin{eqnarray}
\nonumber
u(\V{x},t) & = & \max \left\{ \frac{1}{\kappa^d} \left( 1 - \frac{\tilde{\V{x}}^T \tilde{\V{x}}}{r_0^2 \kappa^2} \right)^{\frac{1}{m}}, \quad 0 \right\} \\
\label{soln-apme}
& = & \max \left\{ \frac{1}{\kappa^d} \left( 1 - \frac{\V{x}^T \D^{-1} \V{x}}{r_0^2 \kappa^2} \right)^{\frac{1}{m}}, \quad 0 \right\}
\end{eqnarray}
where $\kappa$, $\beta$, and $t_0$ are the same as those defined in \eqref{soln-para}.

\begin{rem}
	If the initial solution of APME is taken as \eqref{soln-pme-ini} instead of \eqref{soln-apme-ini}, the boundary of the compact support is not a circle in the computational domain $\Omega_c$. Therefore, the evolution of the solution cannot be described by the solution \eqref{soln-apme}. The difference will be demonstrated in Examples \ref{ex1} and \ref{ex2}.
\end{rem}

\section{Linear finite element formulation}
\label{Sec-fem}

In this section, we briefly describe the linear finite element formulation of IBVP (\ref{ibvp-1}). 
Assume that an affine family of simplicial triangulations $\{ \mathcal{T}_h \}$ is given for the physical
domain $\Omega$, and define
\[
U_0 = \{ v \in H^1(\Omega) \; | \: v|_{\partial \Omega} = 0\}.
\]
Denote the linear finite element space associated with mesh $\mathcal{T}_h$ by $U_0^h$. Then a linear finite element solution $u^h(t) \in U_0^h$  for $t \in (0, T]$ to IBVP (\ref{ibvp-1}) is defined by
\begin{equation}
	\label{fem-form}
	\int_{\Omega} \frac{\partial u^h}{\partial t} \, v^h  d\V{x} + 
	\int_{\Omega} (\nabla v^h)^T \left( (u^h)^{m} \mathbb{D} \right) \nabla u^h  d\V{x}  = 0, 
	\quad \forall v^h \in U_0^h, \quad t \in (t_0,T].
\end{equation} 
Denote the numbers of the elements, vertices, and interior vertices of $\mathcal{T}_h$
by $N$, $N_v$, and $N_{vi}$, respectively. Assume that the first $N_{vi}$ vertices are the interior vertices. Then $U_0^h$ and $u^h$ can be expressed as
\bey
&& U_0^h(t) = \text{span} \{ \phi_1, \cdots, \phi_{N_{vi}} \},
\nn \\
\label{soln-approx}
&& u^h = \sum_{j=1}^{N_{vi}} u_j(t) \phi_j + \sum_{j=N_{vi}+1}^{N_{v}} u_j(t) \phi_j,
\eey
where $\phi_j$ is the linear basis function associated with the $j^{\text{th}}$ vertex, $\V{x}_j$, at time $t$. The boundary and initial conditions in (\ref{ibvp-1}) are approximated as
\beq
u_j(t) = 0, \quad j = N_{vi}+1, ..., N_v, 
\label{fem-bc}
\eeq
and
\beq
u_j(t_0) = u_0(\V{x}_j), \quad j = 1, ..., N_v .
\label{fem-ic}
\eeq

Substituting (\ref{soln-approx}) into (\ref{fem-form}), taking $v^h = \phi_i$ ($i=1, ..., N_{vi}$),
and combining the resulting equations with (\ref{fem-bc}), we obtain the linear algebraic system
\beq
\label{fem-sys}
M \, \frac{d \V{u}} {d t} + A(u^h) \, \V{u} = \V{0},
\eeq
where $\V{u} = (u_1,..., u_{N_{vi}}, u_{N_{vi}+1},..., u_{N_v})^T$ is the unknown vector and $M$ and $A$
are the mass and stiffness matrices, respectively. The entries of the matrices are given as follows. 
For  $j=1, ..., N_{v}$,
\begin{align}
\label{matM}
& m_{ij} = \begin{cases}
	\int_{\Omega} \phi_j \phi_i \, d\V{x} = \sum\limits_{K \in \Th} \int_K \phi_j \phi_i \, d\V{x}, & i=1, ..., N_{vi} \\
	0, & i=N_{vi}+1, ..., N_{v} 
\end{cases} 
\\
\label{matA}
& a_{ij} = \begin{cases}
	\int_{\Omega} (\nabla \phi_i)^T \; (u^h)^{m} \D \nabla \phi_j \, d\V{x}
	= \sum\limits_{K \in \Th} \int_K (\nabla \phi_i )^T \; (u^h)^{m} \D  \nabla \phi_j \, d\V{x}, & i=1, ..., N_{vi} \\
	\delta_{ij}, \quad i=N_{vi}+1, ..., N_{v}. &
\end{cases}
\end{align}

The system \eqref{fem-sys} is solved using the fifth-order Radau IIA method with a two-step error estimator \cite{GMP04}. The relative and absolute tolerances are chosen as $10^{-6}$ and $10^{-8}$, respectively.

Next, we apply the anisotropic mesh adaptation method to generate the mesh $\Th^n$ for time $t=t^n$ ($n=0,1,\cdots$).

\section{Anisotropic mesh adaptation}
\label{Sec-adap}

In this section we briefly introduce the anisotropic mesh adaptation method used in the computations for this paper. More details about the method can be found in \cite{Hua05b, Hua06, LH10, HR11}. 
We take the so-called $\M$-uniform mesh approach where an adaptive mesh is viewed as a uniform one in the metric specified
by a tensor $\M = \M(\V{x})$ that is assumed to be symmetric and uniformly positive
definite on $\Omega$. It is shown in \cite{Hua06} that an $\M$-uniform mesh $\Th$ satisfies
\begin{align}
\label{eq-1}
&|K| \det(\MK)^{\frac{1}{2}}   =  \frac{\sigma_h}{N}, \quad \forall K \in \Th \\
\label{ali-1}
&\frac{1}{d} \mbox{tr} \left ( (F_K')^T \MK F_K' \right )  = 
\mbox{det} \left ( (F_K')^T \MK F_K' \right )^{\frac{1}{d}}, \quad \forall K \in \Th
\end{align}
where 
\beq
\MK = \frac{1}{|K|} \int_K \M(\V{x}) \, d \V{x},\quad
\sigma_h = \sum_{K \in \Th} \det(\MK)^{\frac{1}{2}} |K|.
\eeq
Condition (\ref{eq-1}) is called as {\em the equidistribution condition} which determines the size of $K$ from $\det(\MK)^{\frac{1}{2}}$,
while (\ref{ali-1}) is called {\em the alignment condition} which controls the shape and orientation of $K$. 

In this paper, we consider three choices of metric tensors for adaptive meshes. The first choice is based on
minimization of the $H^1$ semi-norm of linear interpolation error and is given in \cite{Hua05b} by
\beq
\label{M-adap}
\M_{adap}(K) = \Big \| I + \frac{1}{\alpha_h} | H_K(u^h) | \Big \| ^{\frac{2}{5}}
\det\left( I+\frac{1}{\alpha _{h}}|H_K(u^h)|\right) ^{-\frac{1}{5}}
\left[ I+\frac{1}{\alpha _{h}}|H_K(u^h)|\right] ,
\eeq
where $u^h$ is the finite element solution, $H_K(u^h)$ is a recovered Hessian of $u^h$ over $K$,
$|H_K(u^h)|$ is the eigen-decomposition of $H_K(u^h)$ with the eigenvalues being replaced
by their absolute values, and $\alpha_h$ is a positive regularization parameter.

The second and third choices are related to the diffusion matrix and are defined in \cite{LH10} as
\beq
\M_{DMP}(K) = \DK^{-1}, \quad \forall K \in \Th
\label{M-DMP}
\eeq
and
\beq
\label{M-DMP+adap}
\M_{DMP+adap}(K) = \left ( 1+ \frac{1}{\alpha_h}  B_K \right )^{\frac{2}{d+2}}
\mbox{det} \left ( \DK \right )^{\frac{1}{d}} \DK^{-1},
\eeq
where
\begin{align*}
& \DK = \frac{1}{|K|} \int_K \D(\V{x}) \, d\V{x}, \\
& B_K = \mbox{det} \left ( \DK \right )^{-\frac{1}{d}}
 \|\DK^{-1} \|\cdot   \|\DK |H_K(u^h)| \|^2, \\
& \alpha_h = \left ( \frac{1}{|\Omega|} \sum_{K \in \Th} |K| B_K^{\frac{d}{d+2}} \right )^{\frac{d+2}{d}}.
\end{align*}
The $\M$-uniform meshes associated with $\M_{DMP}$ and $\M_{DMP+adap}$ satisfy the alignment condition \eqref{ali-1}, and the mesh elements are aligned along the principal diffusion direction that is the direction of the eigenvector corresponding to the largest eigenvalue of the diffusion matrix $\D$.

If only fixed mesh is used in the computations, a very fine initial mesh is needed in order to capture the sharp change of solution at the initial free boundary. Even for moving mesh methods, an initial fine mesh is usually needed since the connectivity of the mesh elements is fixed. For adaptive meshes, however, they only need to concentrate elements around the initial free boundary or moving boundary. Therefore, an anisotropic initial mesh will significantly increase the efficiency and accuracy of the computation. The anisotropic initial mesh can be generated by adapting a uniform Delaunay mesh according to a metric tensor $\M$. The adaptation can iterate a few times to generate an initial mesh with good quality. 

\begin{figure}[!thb]
\centering
\includegraphics[width=4.5in]{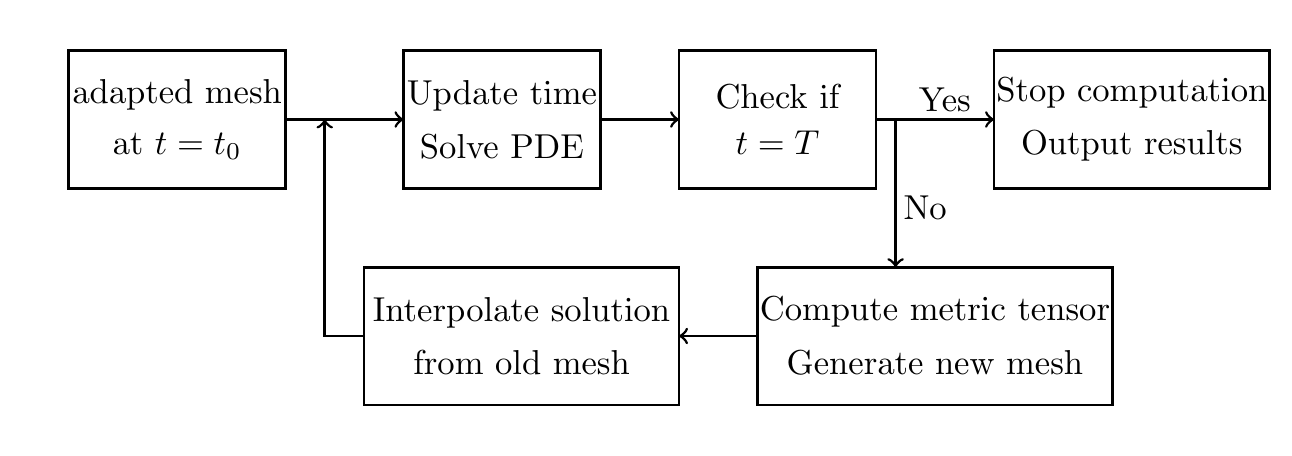}
\caption{Procedures for anisotropic mesh adaptation based on $\M$-uniform mesh approach}
\label{meshadap}
\end{figure}

Starting with the adapted initial mesh, the procedures for anisotropic mesh adaptation based on $\M$-uniform mesh is shown in Fig. \ref{meshadap}. Some numerical results are presented in the next section.

\begin{rem}
	In the mesh adaptation procedures shown in Fig. \ref{meshadap}, we use linear finite element interpolation to map the solution from old mesh to the new adaptive mesh for the current time $t^n$. Then the new mesh is used to solve PDE for the next time $t^{n+1}$. Higher order interpolations can also be used if needed.
\end{rem}

\begin{rem}
	The adaptation for the mesh at a given time $t^n$ can be iterated a few times. However, the solution also needs to be interpolated correspondingly, which may reduce the accuracy of the solution. Therefore, the mesh is only adapted once for each time step in our computations.  
\end{rem}

\begin{rem}
	Comparing to uniform meshes, the condition numbers of mass matrix $M$ in \eqref{matM} and stiffness matrix $A$ in \eqref{matA} for anisotropic meshes are affected by mesh non-uniformity. However, the condition numbers are still bounded. In anisotropic meshes, elements are usually concentrated in a small portion of the physical domain. With Jacobi preconditioning, the condition numbers of $M$ and $A$ for $\M$-uniform meshes are comparable with Delaunay meshes \cite{KHX14} and the impact on iterative convergence is not significant.
\end{rem}

\section{Numerical results}
\label{Sec-results}
In this section we present numerical results obtained in two dimensions for three examples to demonstrate the different behavior between PME and APME as well as the significance of mesh adaptation in the computations. For comparison purpose, we consider four types of meshes. One is a fixed mesh that does not change through out the computations. The other three meshes are $\M_{adap}$ mesh, $\M_{DMP}$ mesh, and $\M_{DMP+adap}$ mesh that have been introduced in Section \ref{Sec-adap}. 

The fixed mesh is generated by splitting the rectangle domain into sub-rectangles, then each sub-rectangles are divided into four triangles by the diagonal lines. The adaptive meshes are generated using the Bidimensional Anisotropic Mesh Generator (BAMG) developed by Hecht \cite{bamg} based on the Delaunay triangulation and local node movement.

\begin{exam}
\label{ex1}
The first example is in the form of IBVP (\ref{ibvp-1}) for $m=1$ on domain $\Omega = [-3,3]^2$ with diffusion
$\D = \begin{pmatrix} 5.5 & 4.5 \\ 4.5& 5.5 \end{pmatrix}$ and $u_0$ given in \eqref{soln-apme-ini} with $r_0=0.5$. The initial time is $t_0=0.03125$ by \eqref{soln-para}.   
The diffusion matrix has eigenvalues $10$ and $1$, and the principal eigenvectors are in the northeast direction. The boundary of the support region for the initial solution is an ellipse in the physical domain $\Omega$. 

Fig. \ref{ex1-mesh-ini} shows the four different initial meshes obtained using different metric tensors for the computations. For $\M_{adap}$ mesh, if not specified otherwise, the regularization parameter in \eqref{M-adap} is chosen as $\alpha_h=0.01$. As can be seen, the elements in both $\M_{DMP}$ mesh and $\M_{DMP+adap}$ mesh are aligned along the north-east direction. 
\begin{figure}[!thb]
\centering
\hbox{
\begin{minipage}{2.5in}
\includegraphics[width=2.5in]{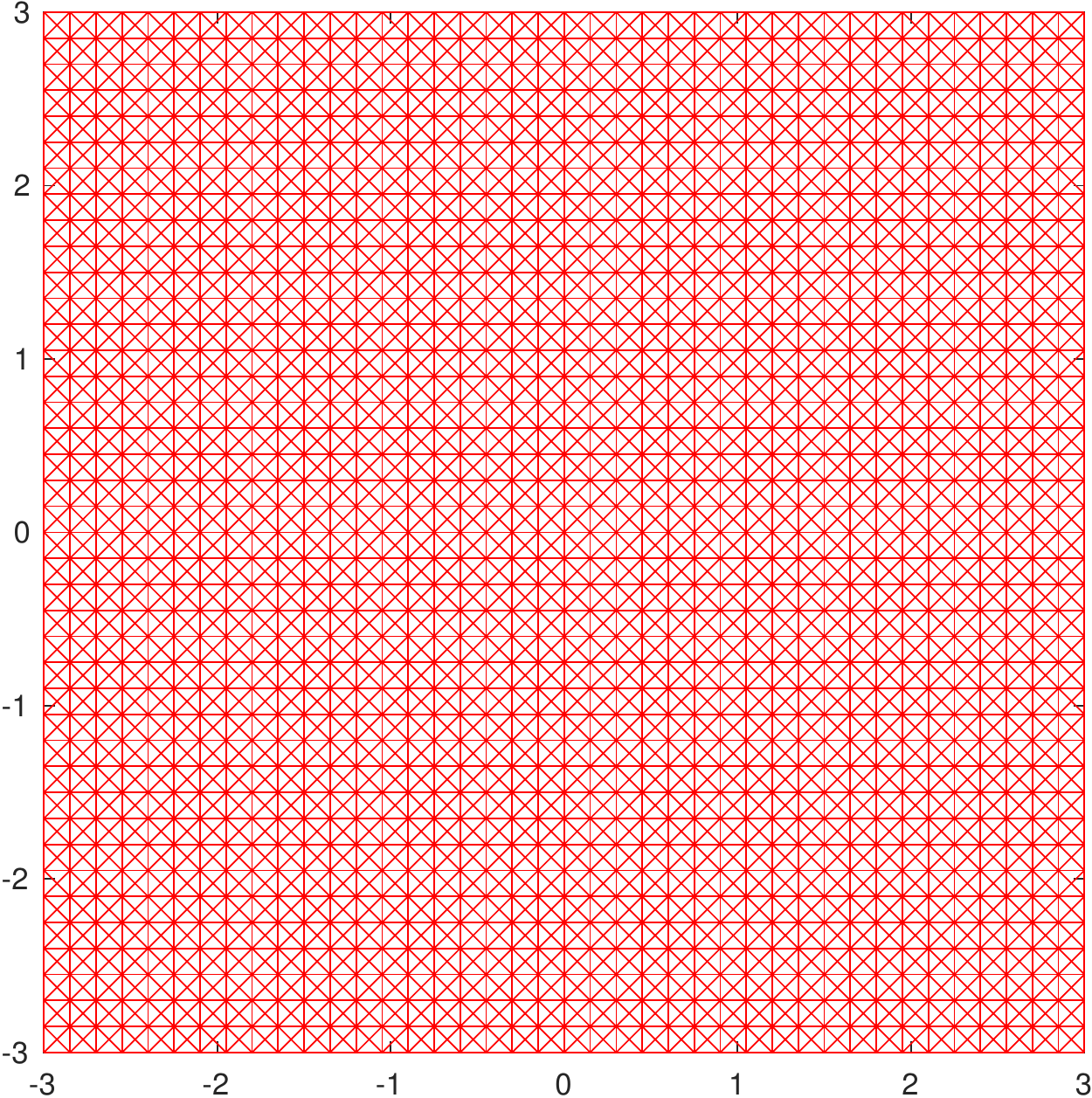}
\centerline{(a): Fixed mesh, $N=6,400$}
\end{minipage}
\begin{minipage}{2.5in}
\includegraphics[width=2.5in]{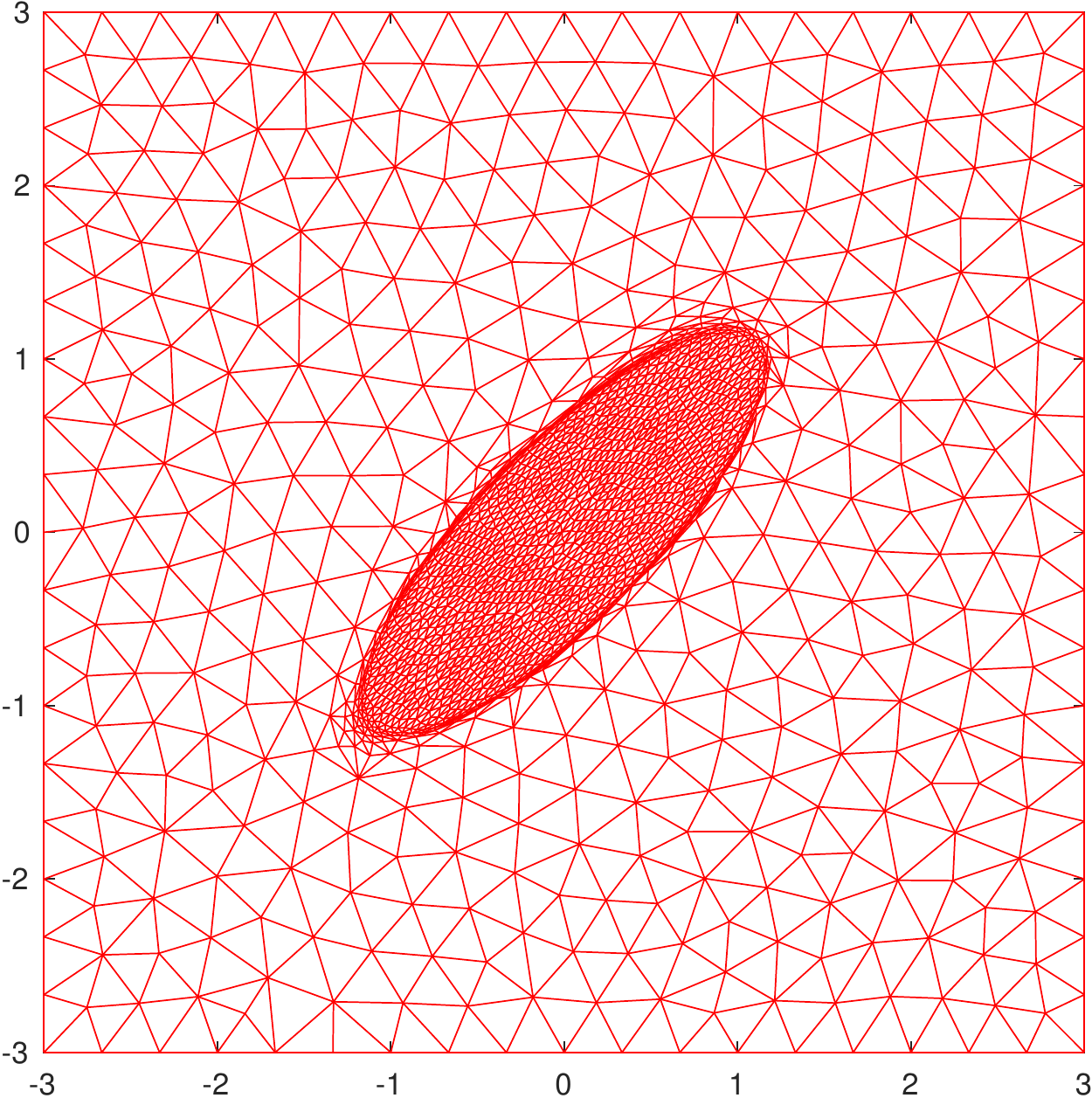}
\centerline{(b): $\M_{adap}$ mesh, $N=3,260$}
\end{minipage}
}
\vspace{5mm}
\hbox{
\begin{minipage}{2.5in}
\includegraphics[width=2.5in]{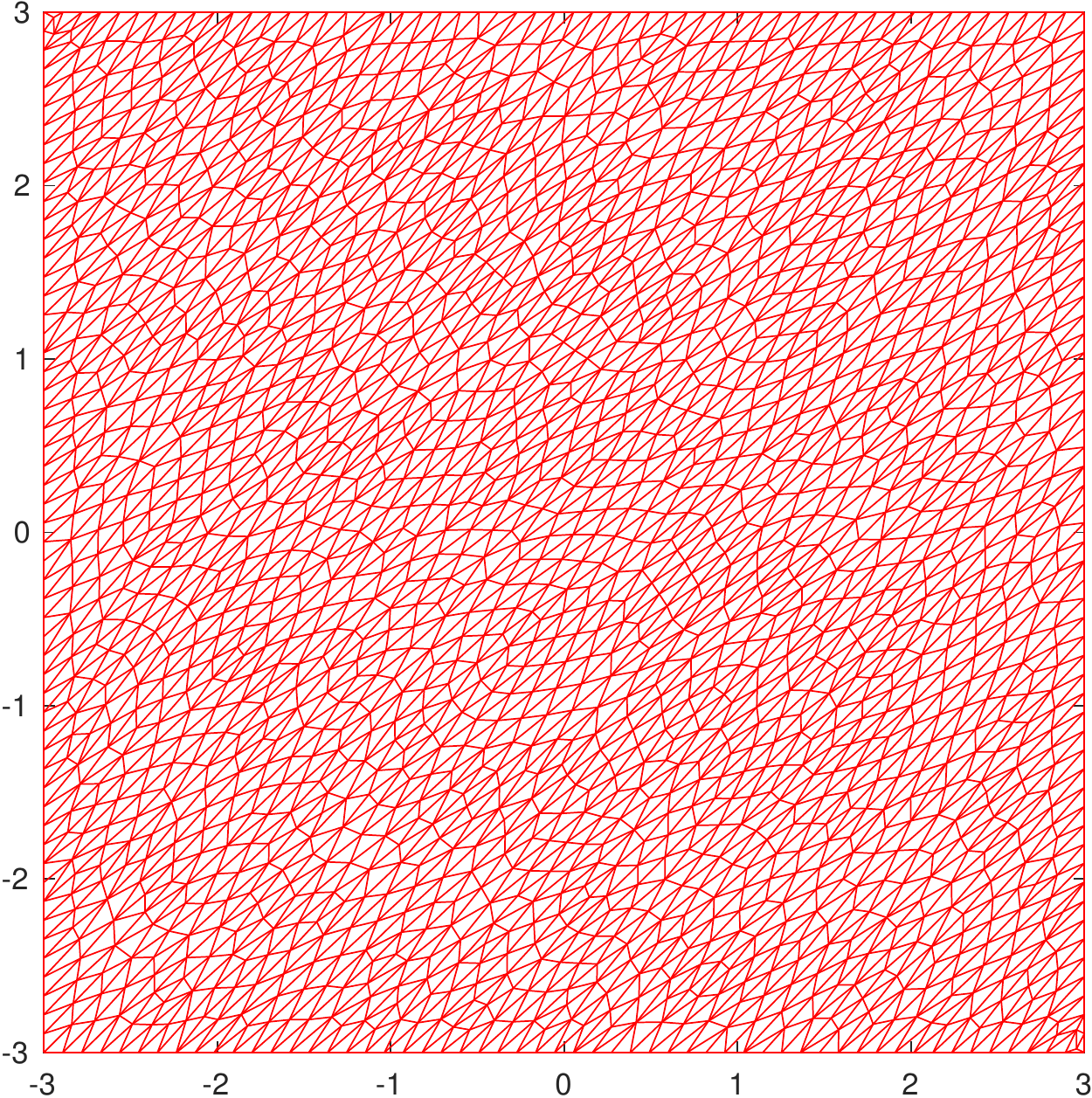}
\centerline{(c):  $\M_{DMP}$ mesh, $N=3,846$}
\end{minipage}
\begin{minipage}{2.5in}
\includegraphics[width=2.5in]{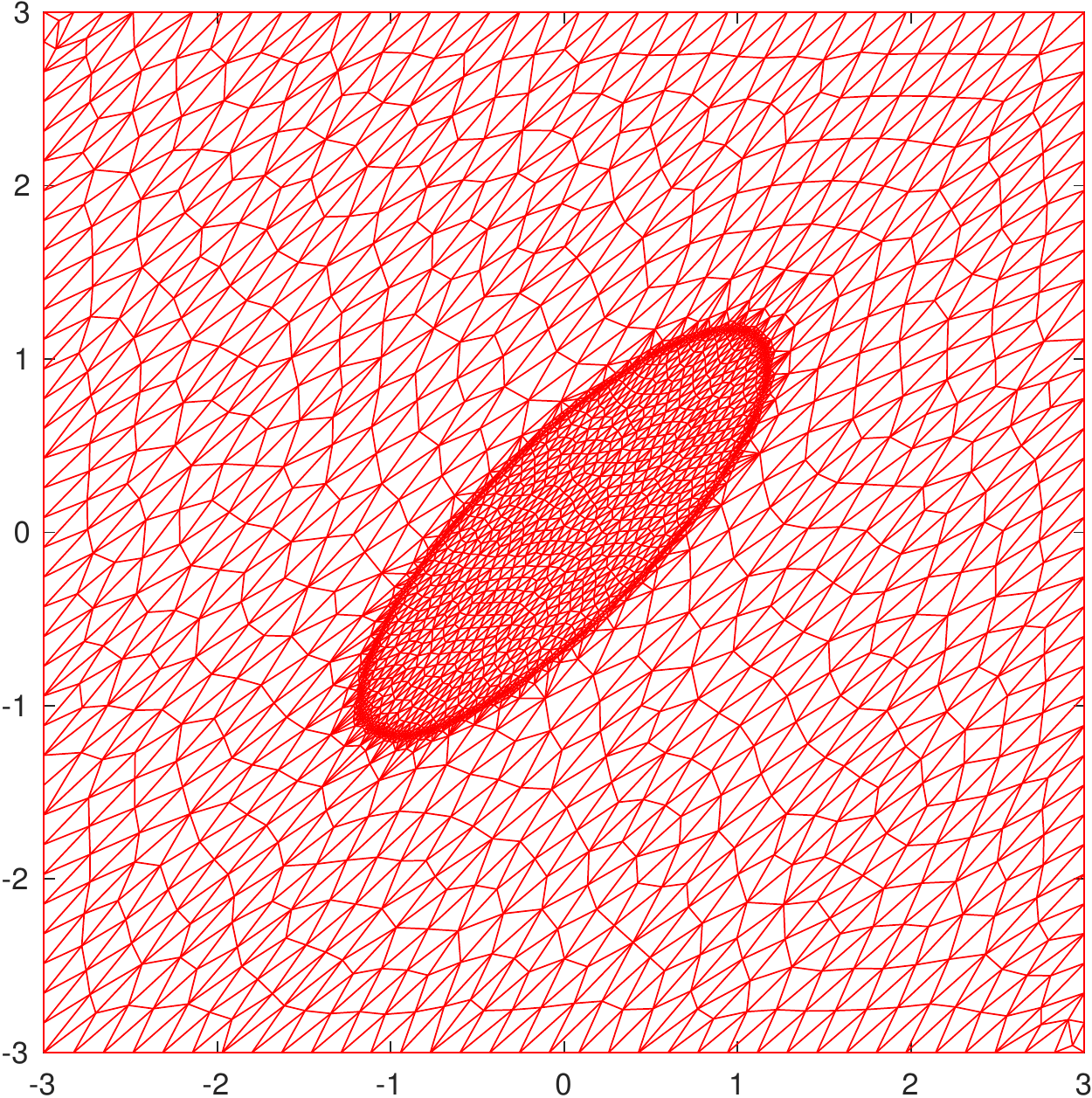}
\centerline{(d): $\M_{DMP+adap}$ mesh, $N=4,522$}
\end{minipage}
}
\caption{Example \ref{ex1}. Initial meshes generated using different metric tensors.}
\label{ex1-mesh-ini}
\end{figure}

The initial solution at $t=t_0$ and the final solution at $t=0.2$ are shown in Fig. \ref{ex1-soln}. The $\M_{adap}$ and $\M_{DMP+adap}$ meshes at different times are shown in Fig. \ref{ex1-mesh}. Comparing with the initial meshes in Fig. \ref{ex1-mesh-ini}(b) and (d), it is clear that the free boundary moves outward in the shape of an ellipse in the physical domain and the meshes are adapted accordingly.   

\begin{figure}[!thb]
\centering
\hbox{
\begin{minipage}{2.5in}
\includegraphics[width=2.5in]{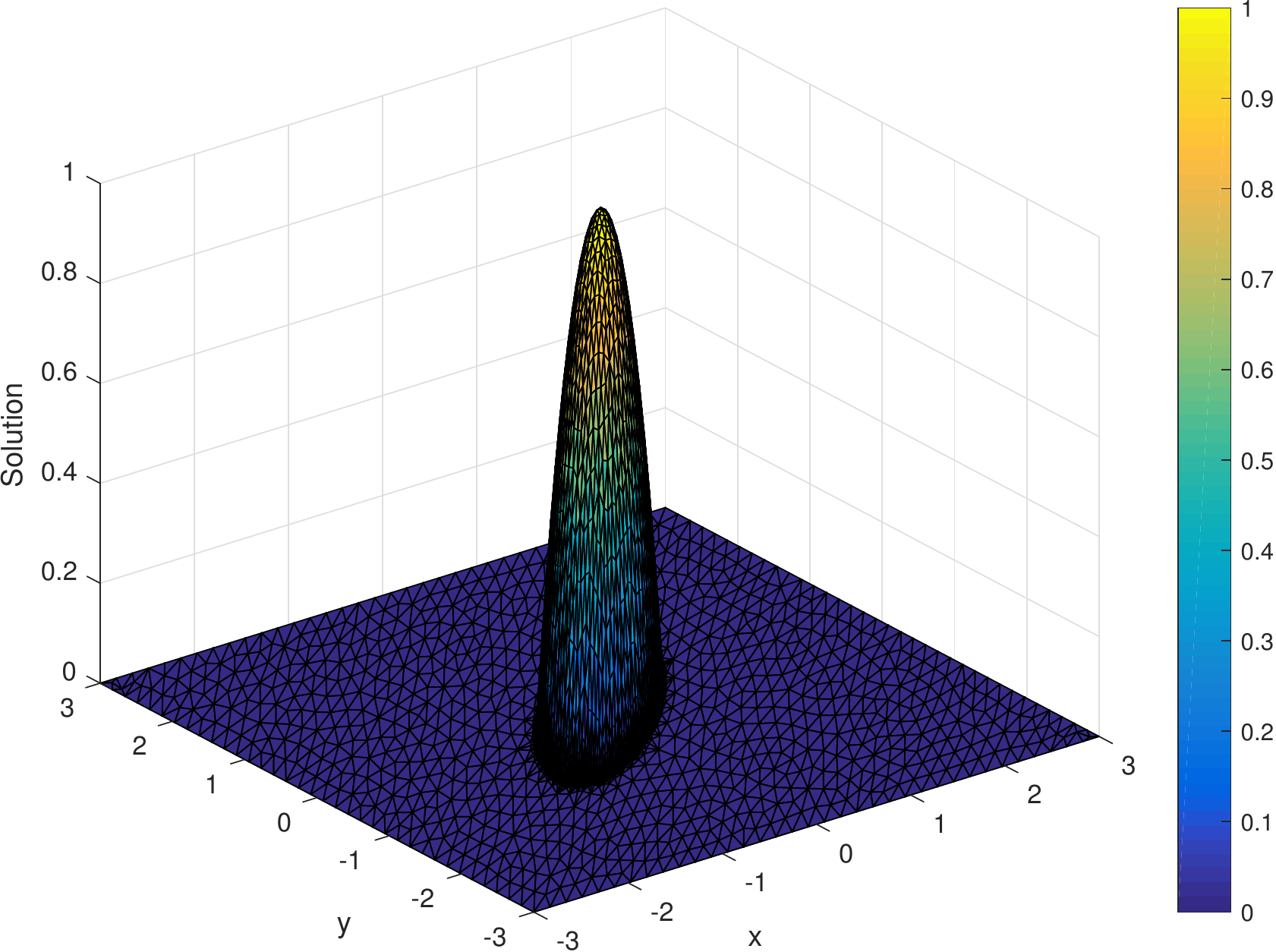}
\centerline{(a): $t=t_0$ }
\end{minipage}
\begin{minipage}{2.5in}
\includegraphics[width=2.5in]{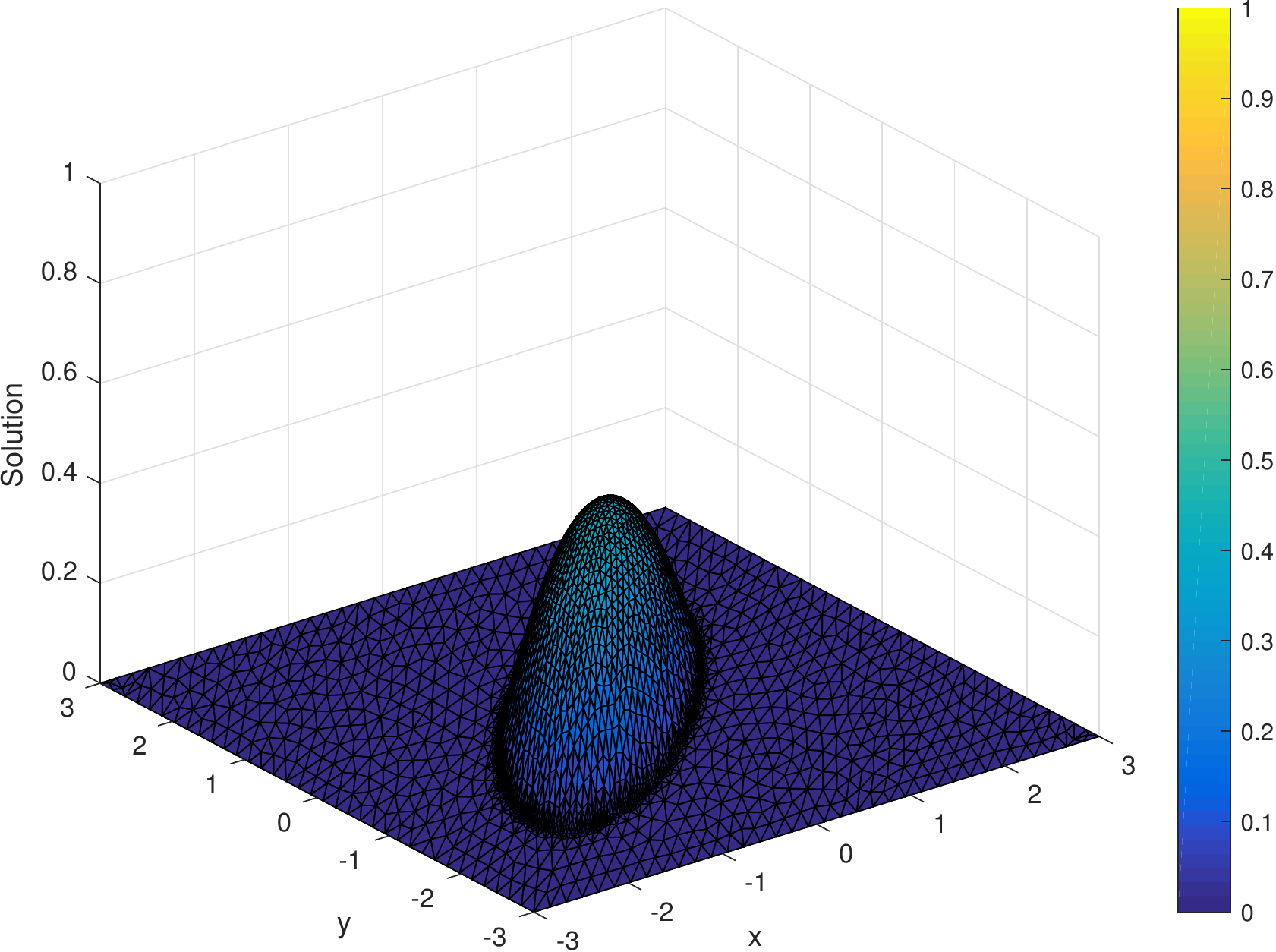}
\centerline{(b): $t=0.2$}
\end{minipage}
}
\caption{Example \ref{ex1}. Numerical solutions at different times.}
\label{ex1-soln}
\end{figure}

\begin{figure}[!thb]
\centering
\hbox{
\begin{minipage}{2.5in}
\includegraphics[width=2.5in]{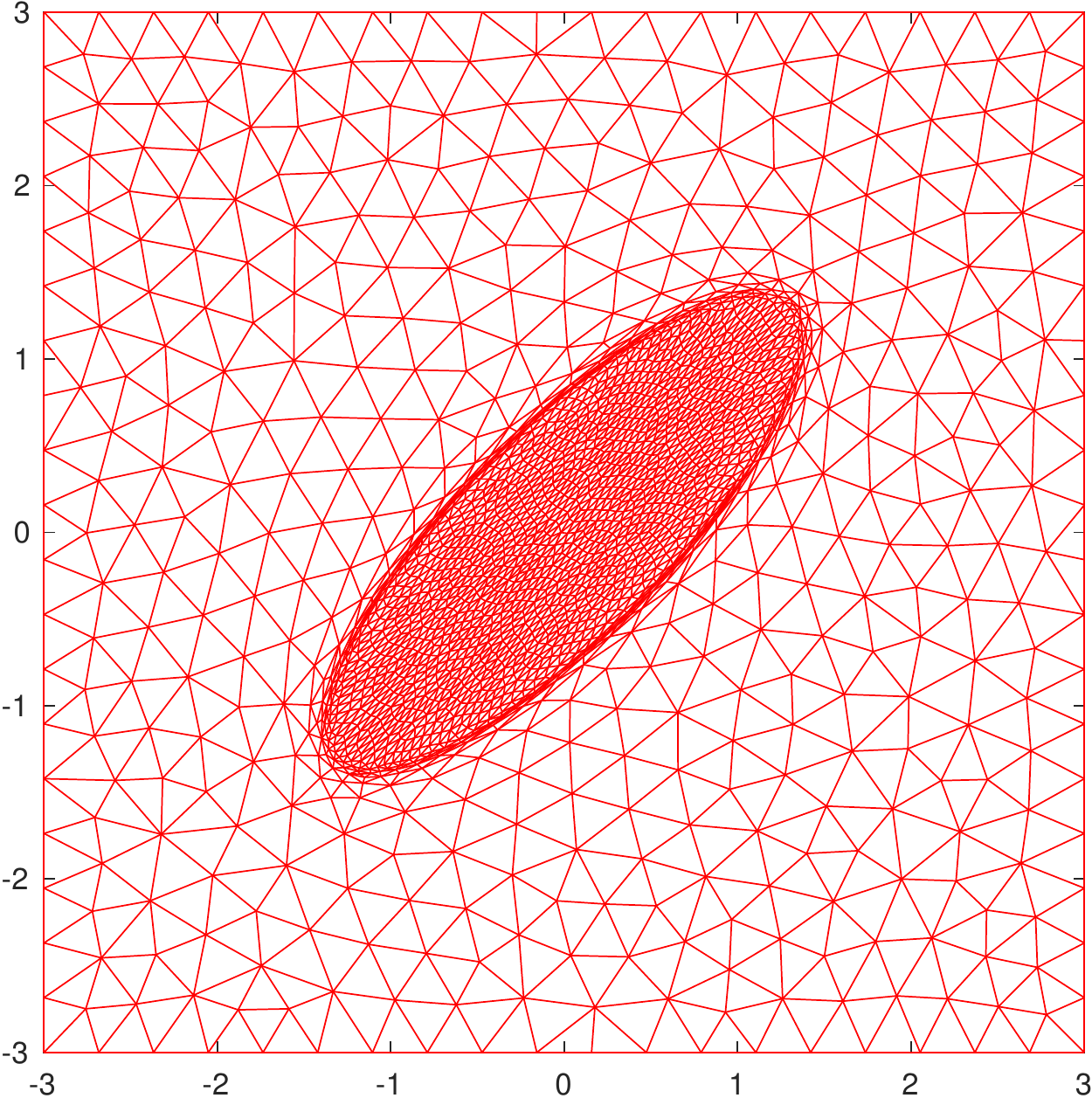}
\centerline{(a): $t=0.06$, $\M_{adap}$ mesh }
\end{minipage}
\begin{minipage}{2.5in}
\includegraphics[width=2.5in]{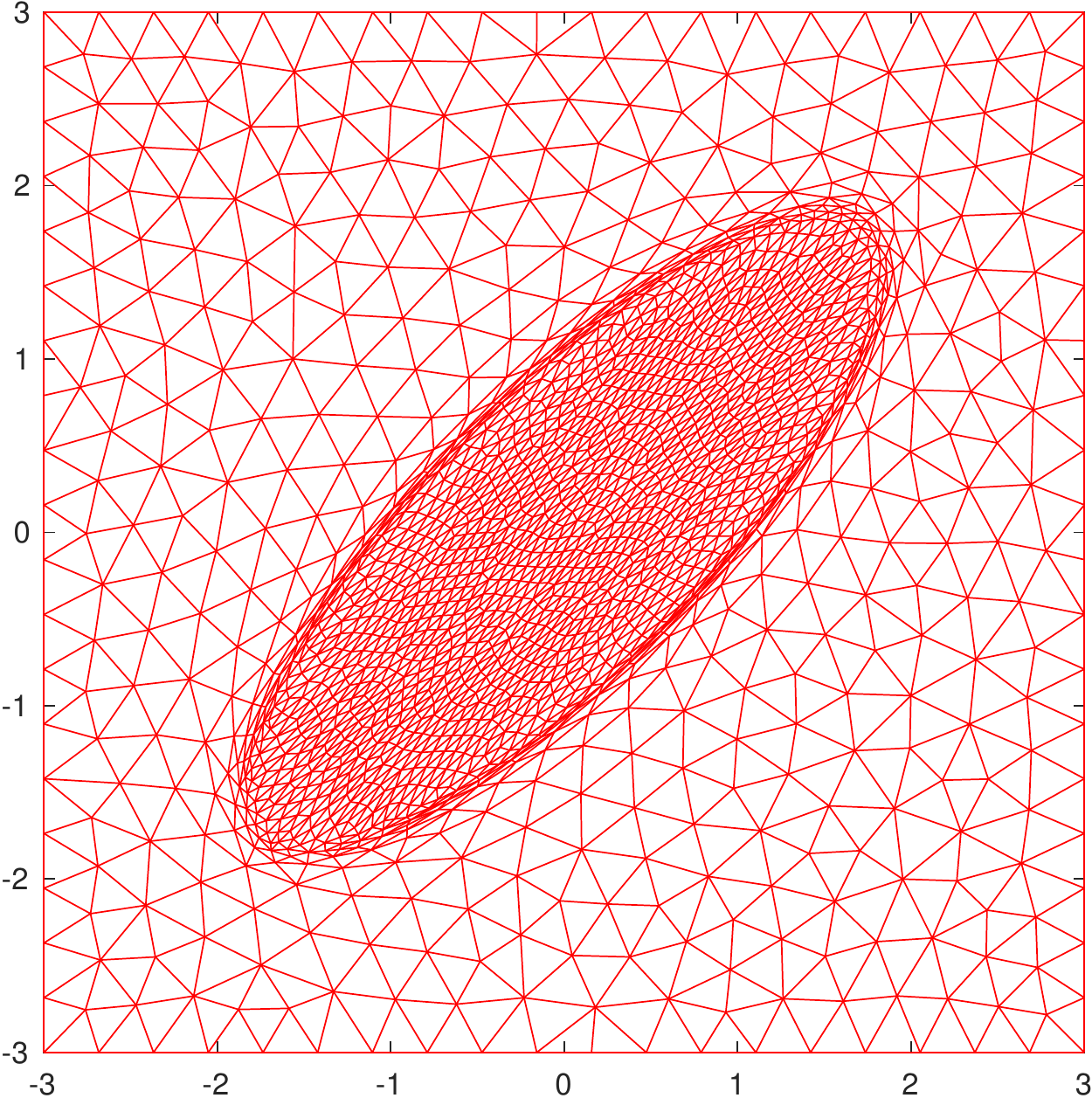}
\centerline{(b): $t=0.2$, $\M_{adap}$ mesh}
\end{minipage}
}
\vspace{5mm}
\hbox{
\begin{minipage}{2.5in}
\includegraphics[width=2.5in]{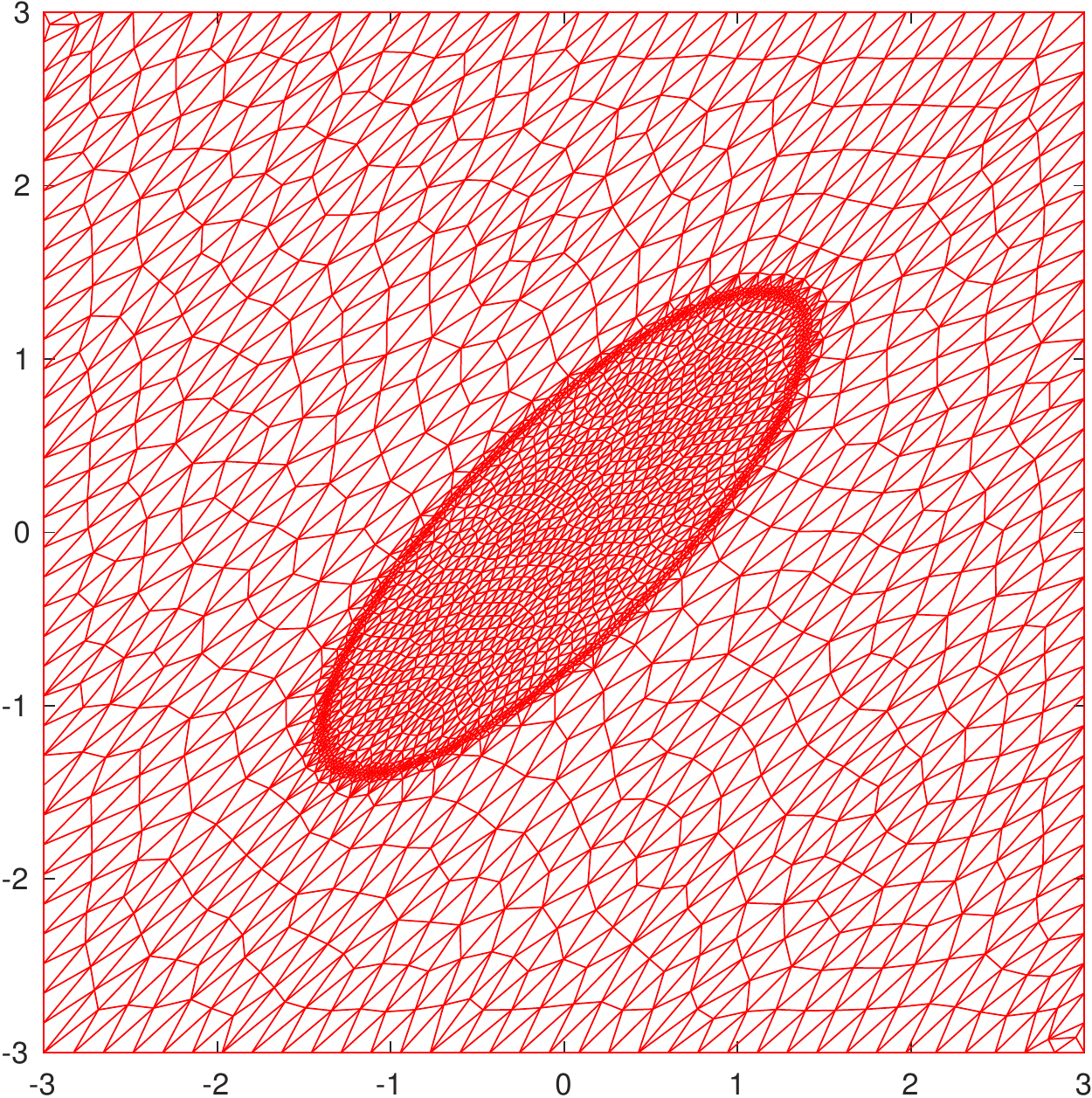}
\centerline{(c): $t=0.06$, $\M_{DMP+adap}$ mesh}
\end{minipage}
\begin{minipage}{2.5in}
\includegraphics[width=2.5in]{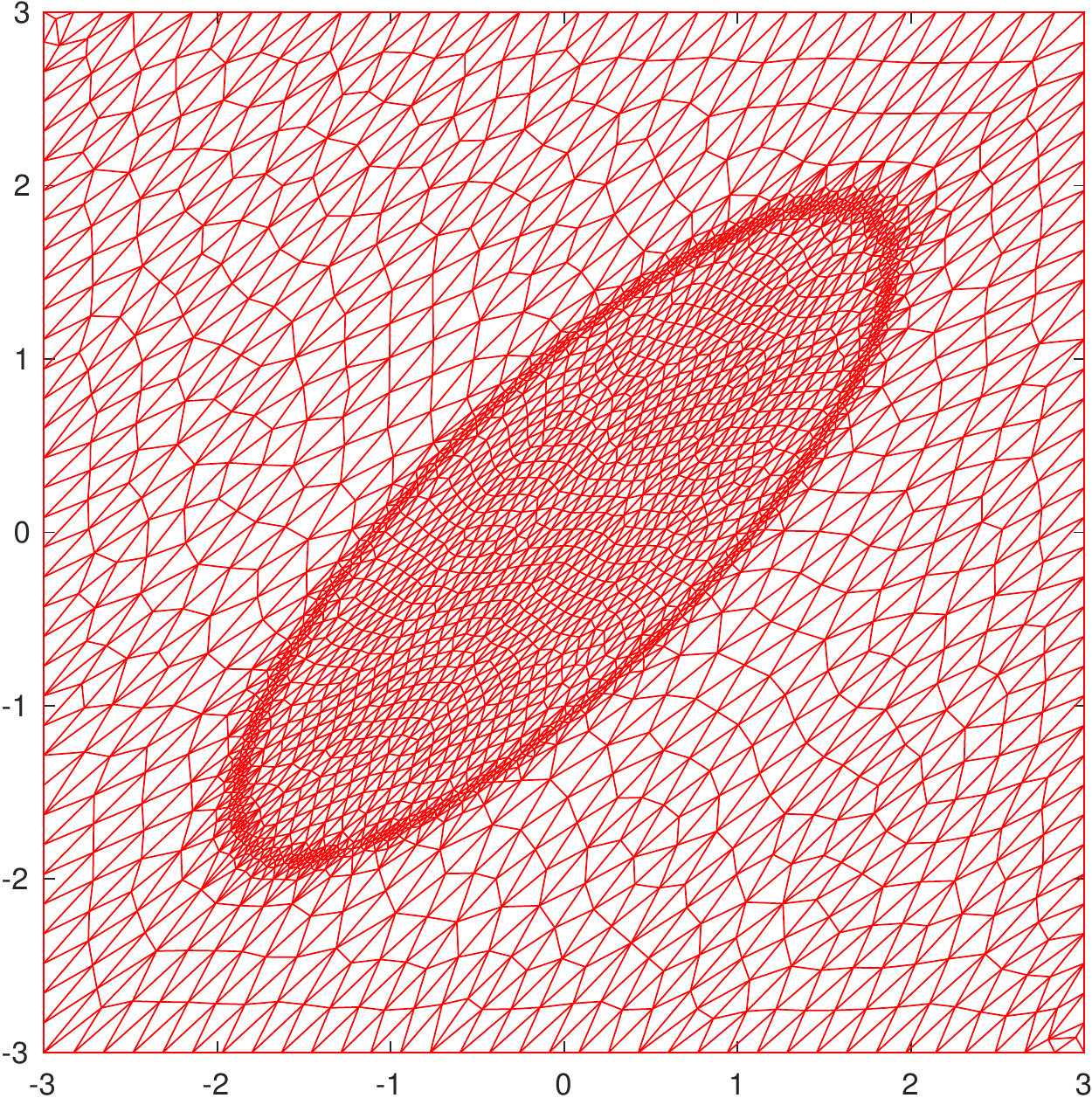}
\centerline{(d): $t=0.2$, $\M_{DMP+adap}$ mesh}
\end{minipage}
}
\caption{Example \ref{ex1}. $\M_{adap}$ and $\M_{DMP+adap}$ meshes at different times.}
\label{ex1-mesh}
\end{figure}

\begin{rem}
From the definition of metric tensors \eqref{M-adap} and \eqref{M-DMP+adap}, the smaller $\alpha_h$ (hence larger $1/\alpha_h$), the more elements will be concentrated around the region with sharp change of solution. Fig. \ref{ex1-meshb} shows $\M_{adap}$ meshes obtained using $\alpha_h=1$ at $t=t_0$ and $t=0.06$. Comparing with Fig. \ref{ex1-mesh-ini}(b) and Fig. \ref{ex1-mesh}(a), it is clear that $\alpha_h=0.01$ provides better adaptation than $\alpha_h=1$ for $\M_{adap}$.
\end{rem}

\begin{figure}[!thb]
\centering
\hbox{
\begin{minipage}{2.5in}
\includegraphics[width=2.5in]{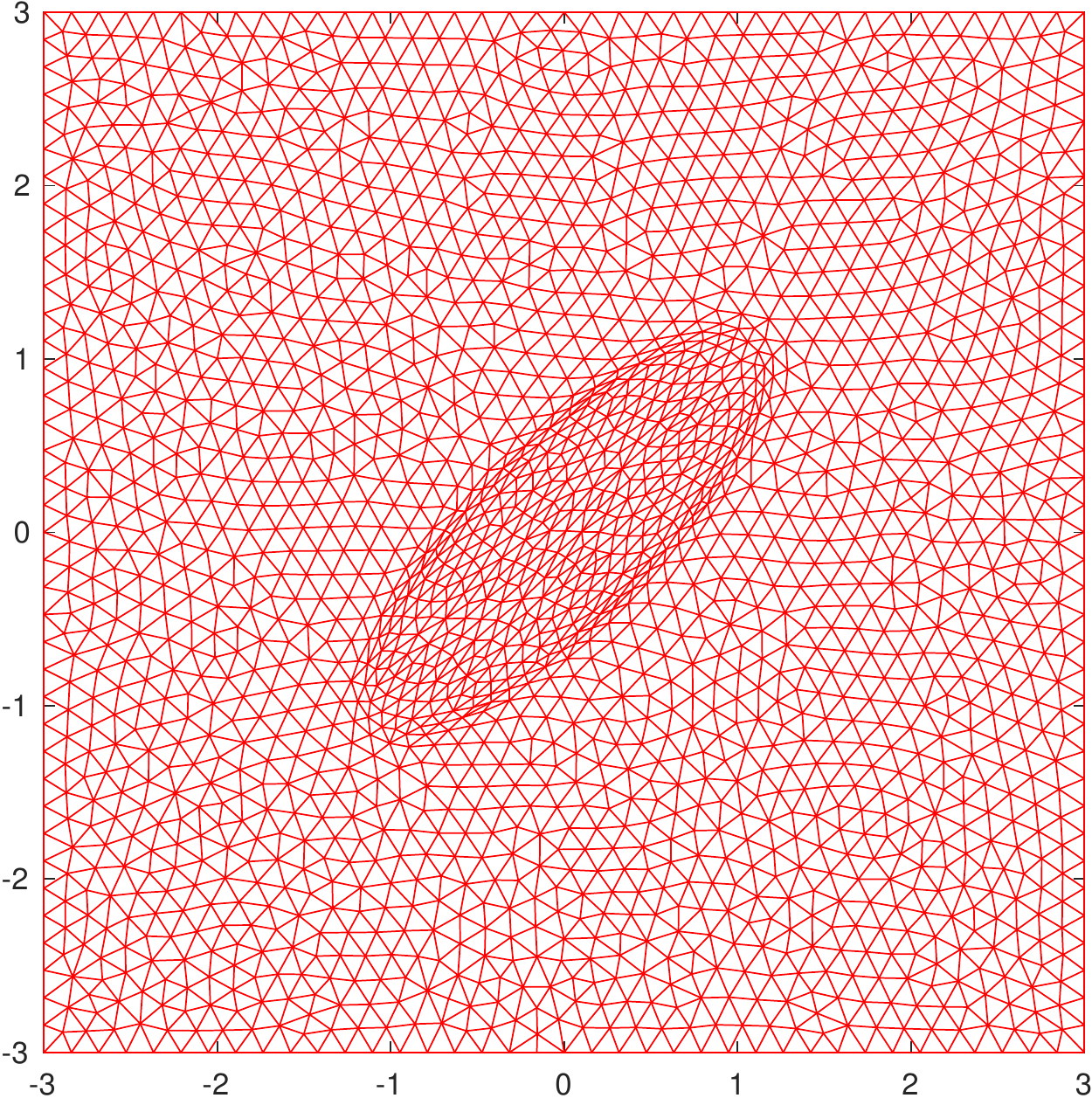}
\centerline{(a): $t=t_0$ }
\end{minipage}
\begin{minipage}{2.5in}
\includegraphics[width=2.5in]{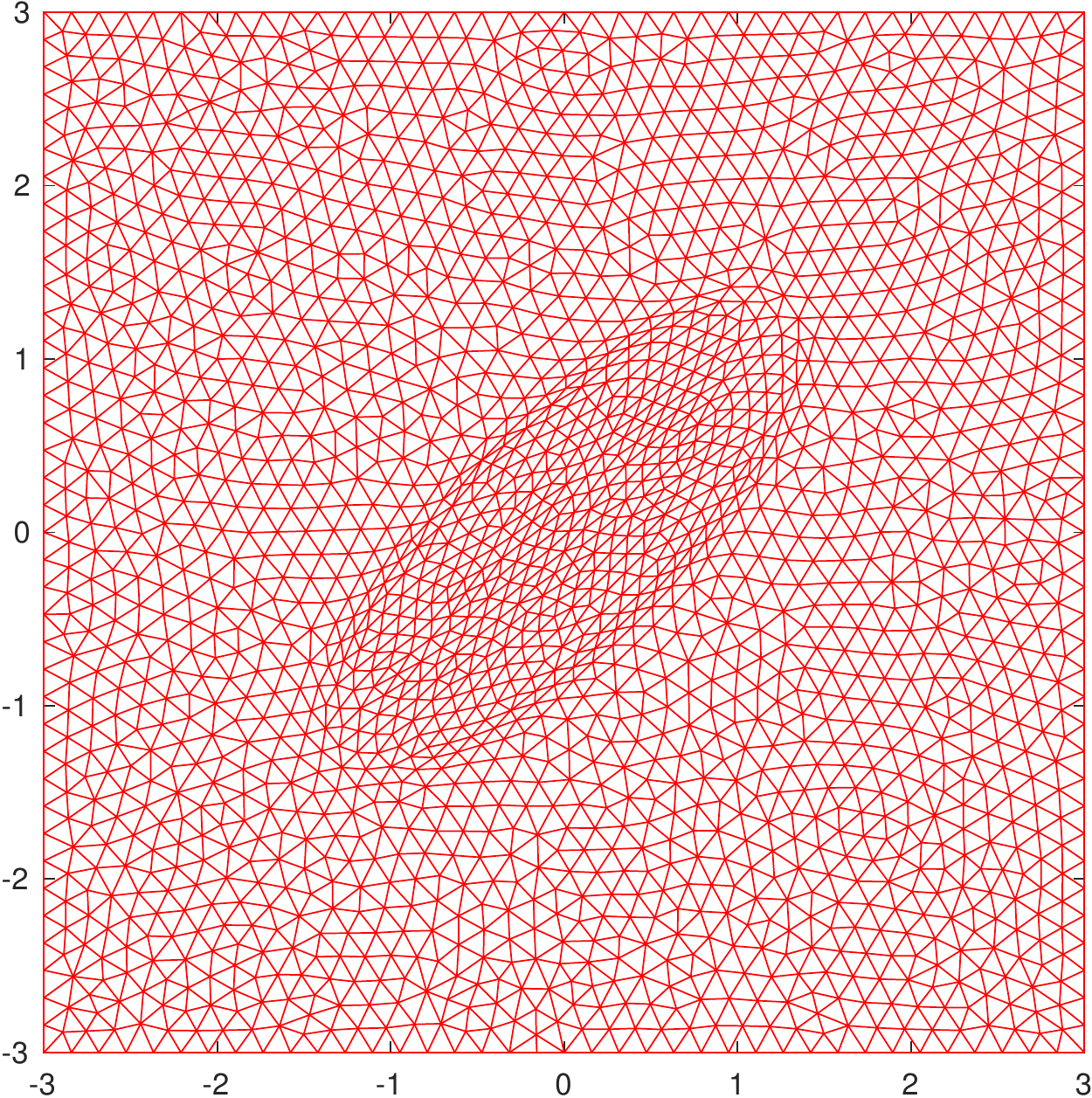}
\centerline{(b): $t=0.06$}
\end{minipage}
}
\caption{Example \ref{ex2}. $\M_{adap}$ meshes obtained with $\alpha_h=1$ at different times.}
\label{ex1-meshb}
\end{figure}

For convergence of solution errors, we choose the final time as $T = (t_0 + 0.1)/2 = 0.065625$. 
The $L^2$ norm of the solution errors obtained from different meshes are plotted in Fig. \ref{ex1-error}, where $\M_{adap,1}$ and $\M_{adap,0.01}$ represent the results from $\M_{adap}$ meshes obtained using $\alpha_h=1$ and $\alpha_h=0.01$, respectively. 
For comparison purpose, the first order and second order reference rates are also plotted in Fig. \ref{ex1-error} and computed as $0.1/\sqrt{N}$ and $1/N$, respectively. 
  
\begin{figure}[!thb]
\centering
\includegraphics[width=4in]{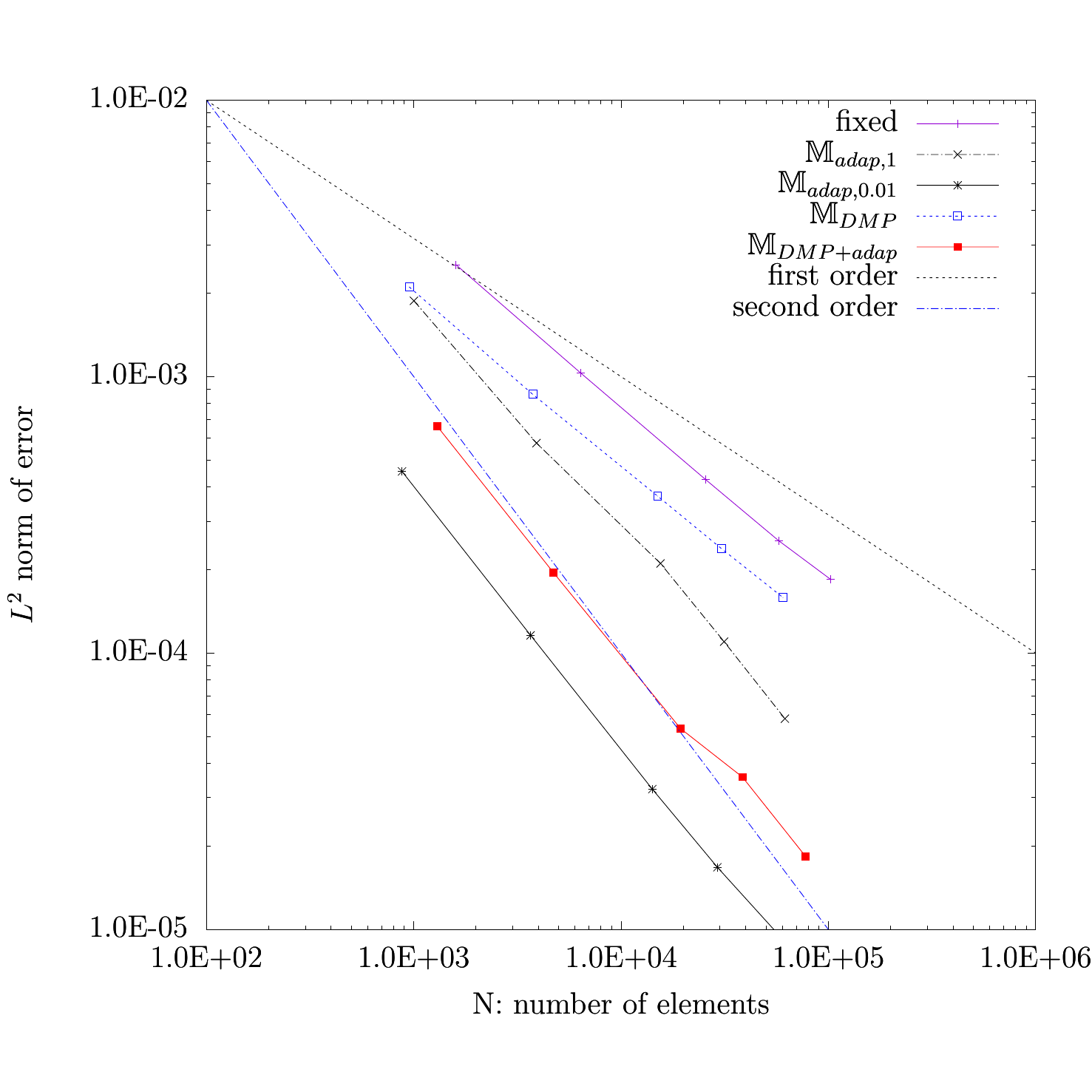}
\caption{Example \ref{ex1}. $L^2$ norm of solution errors obtained from different meshes.}
\label{ex1-error}
\end{figure}

As can be seen from Fig. \ref{ex1-error}, the solution errors obtained from fixed mesh and $\M_{DMP}$ mesh have first order convergence, with $\M_{DMP}$ mesh providing smaller errors than fixed mesh. On the other hand, the errors obtained from $\M_{adap}$ (both $\M_{adap,1}$ and $\M_{adap,0.01}$) and $\M_{DMP+adap}$ meshes have second order convergence, with $\M_{adap,0.01}$ mesh providing smallest errors. 

It is interesting to observe that with slight adaptation using $\M_{adap}$ meshes with $\alpha_h=1$ (that is, $\M_{adap,1}$), the convergence rate of the solution errors has already been improved to second order. With better adaptation using $\alpha_h=0.01$, the solution errors are not only of second order convergence but also much smaller than others. 

The result demonstrates that adaptive meshes help the numerical solutions to achieve higher order convergence. In the mean time, $\M_{DMP+adap}$ mesh also performs well as a combination of $\M_{DMP}$ and $\M_{adap}$, and the errors are between those from $\M_{adap,1}$ and $\M_{adap,0.01}$ meshes.  

\begin{rem}
	Similar to the effect of $\alpha_h$ for $\M_{adap}$ meshes, the regularization parameter $\alpha_h$ in \eqref{M-DMP+adap} for $\M_{DMP+adap}$ meshes can also be adjusted to improve the numerical solutions.  
\end{rem} 

\end{exam}

\vspace{10pt}

\begin{exam}
\label{ex2}
The second example is the same as Example \ref{ex1} except that the initial solution $u_0$ is defined as in \eqref{soln-pme-ini} where the initial free boundary is a circle with radius $r_0$ in the physical domain $\Omega$. The purpose of this example is to show the different behavior between our APME and general PME. For PME, as we know from the exact solution \eqref{soln-pme}, the free boundary will move outward in a shape of a circle. However, it is not the case for APME. 

Fig. \ref{ex2-mesh} shows the initial and final meshes adapted from metric tensors $\M_{adap}$ with $\alpha_h=0.01$ and $\M_{DMP+adap}$. The solutions at different times are shown in Fig. \ref{ex2-soln}. 
As can be seen from the meshes in Fig. \ref{ex2-mesh}, the free boundary varies from a circle gradually to an ellipse. At a given time, the elliptical boundary is different than the ones in Example \ref{ex1} as shown in Fig. \ref{ex1-mesh} (b) and (d). In particular, the eccentricity of the elliptical boundary is smaller than that in Example \ref{ex1}.  The solution at $t=0.2$ shown in Fig. \ref{ex2-soln}(b) is also different than the one in Example \ref{ex1} as shown in Fig. \ref{ex1-soln}(b). 

\begin{figure}[!thb]
\centering
\hbox{
\begin{minipage}{2.5in}
\includegraphics[width=2.5in]{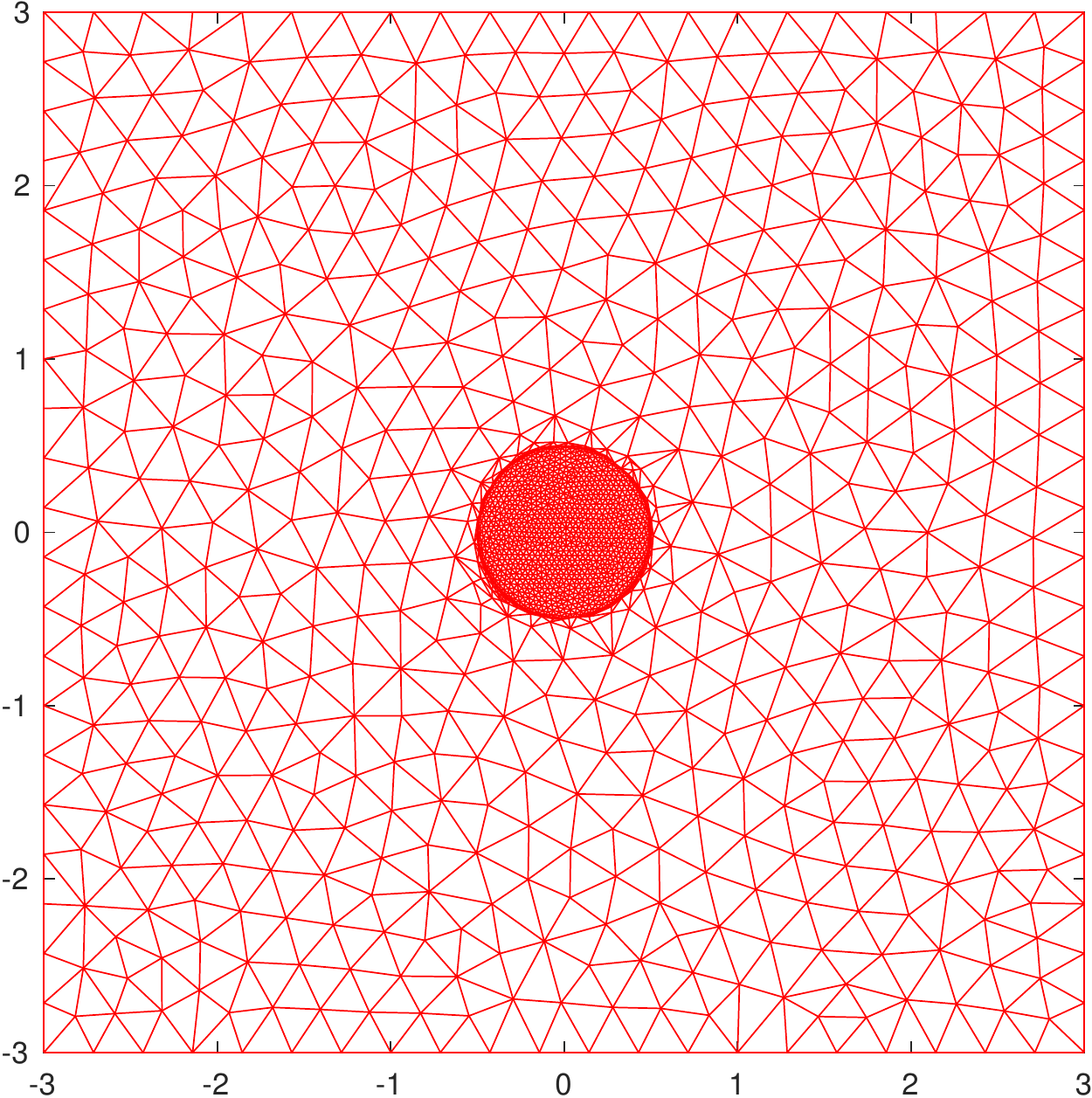}
\centerline{(a): $\M_{adap}$ mesh at $t=t_0$}
\end{minipage}
\begin{minipage}{2.5in}
\includegraphics[width=2.5in]{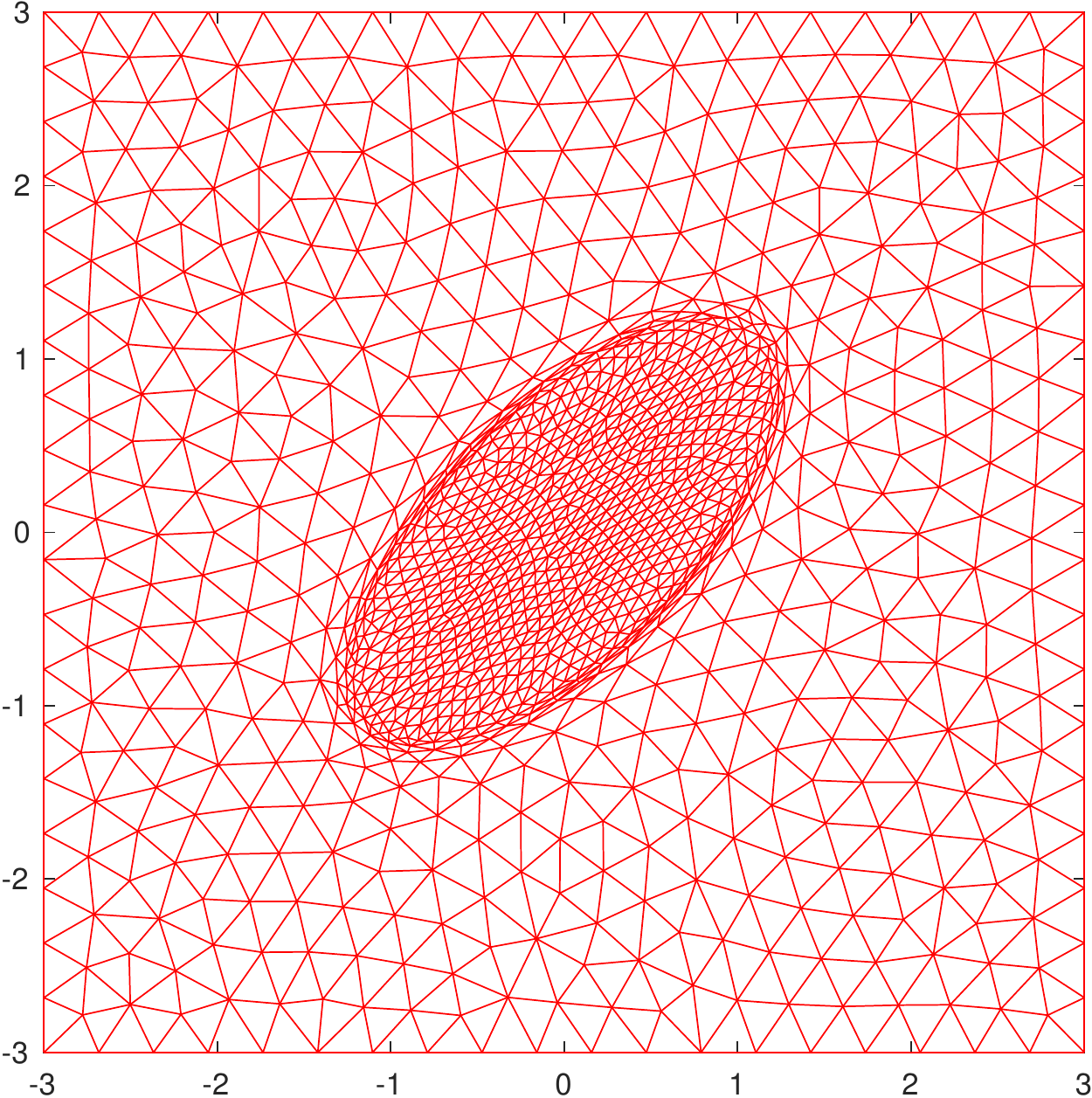}
\centerline{(b): $\M_{adap}$ mesh at $t=0.2$}
\end{minipage}
}
\vspace{5mm}
\hbox{
\begin{minipage}{2.5in}
\includegraphics[width=2.5in]{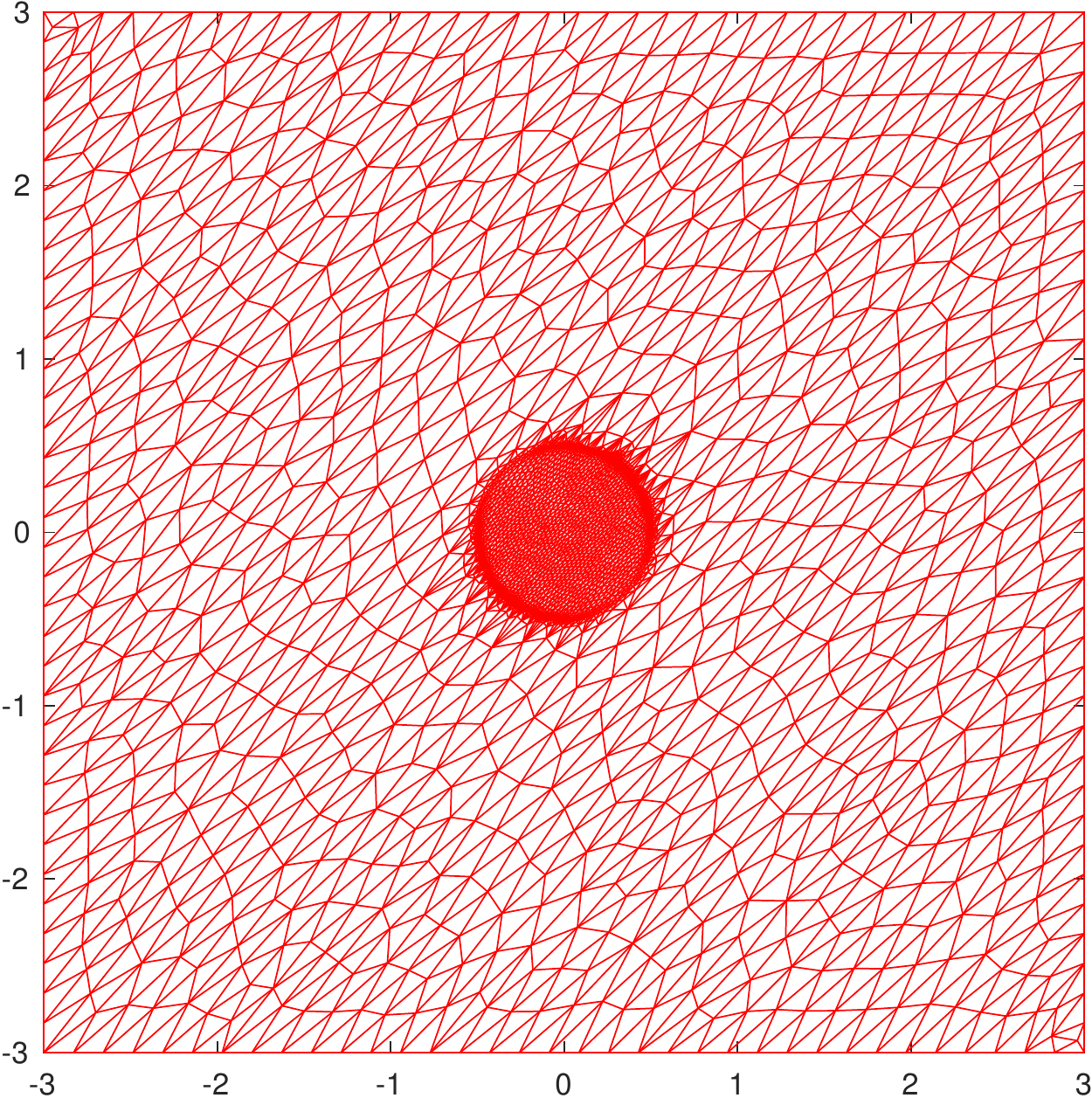}
\centerline{(c):  $\M_{DMP+adap}$ mesh at $t=t_0$}
\end{minipage}
\begin{minipage}{2.5in}
\includegraphics[width=2.5in]{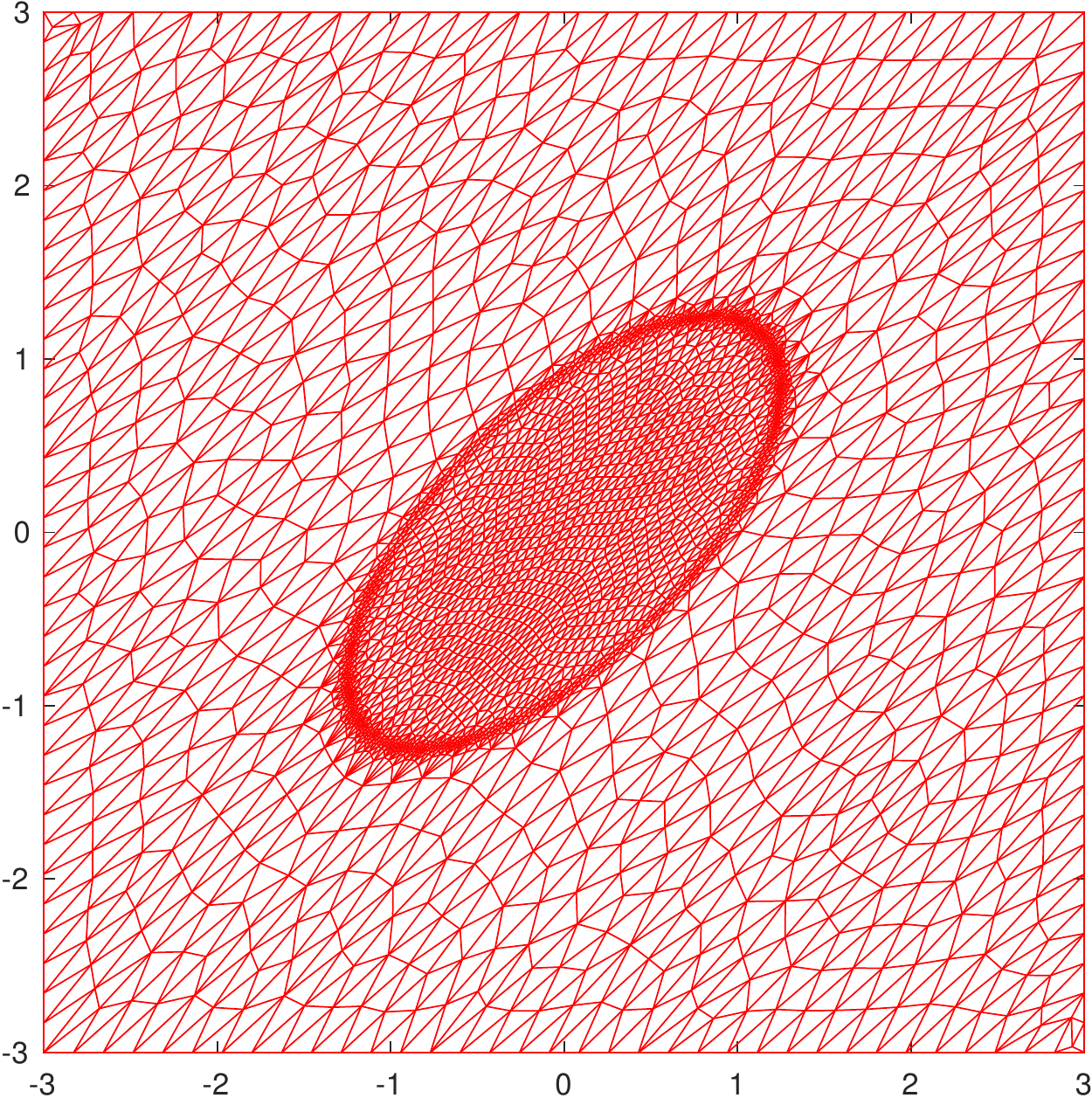}
\centerline{(d): $\M_{DMP+adap}$ at $t=0.2$}
\end{minipage}
}
\caption{Example \ref{ex2}. $\M_{adap}$ and $\M_{DMP+adap}$ meshes at different times, $m=1$.}
\label{ex2-mesh}
\end{figure}

\begin{figure}[!thb]
\centering
\hbox{
\begin{minipage}{2.5in}
\includegraphics[width=2.5in]{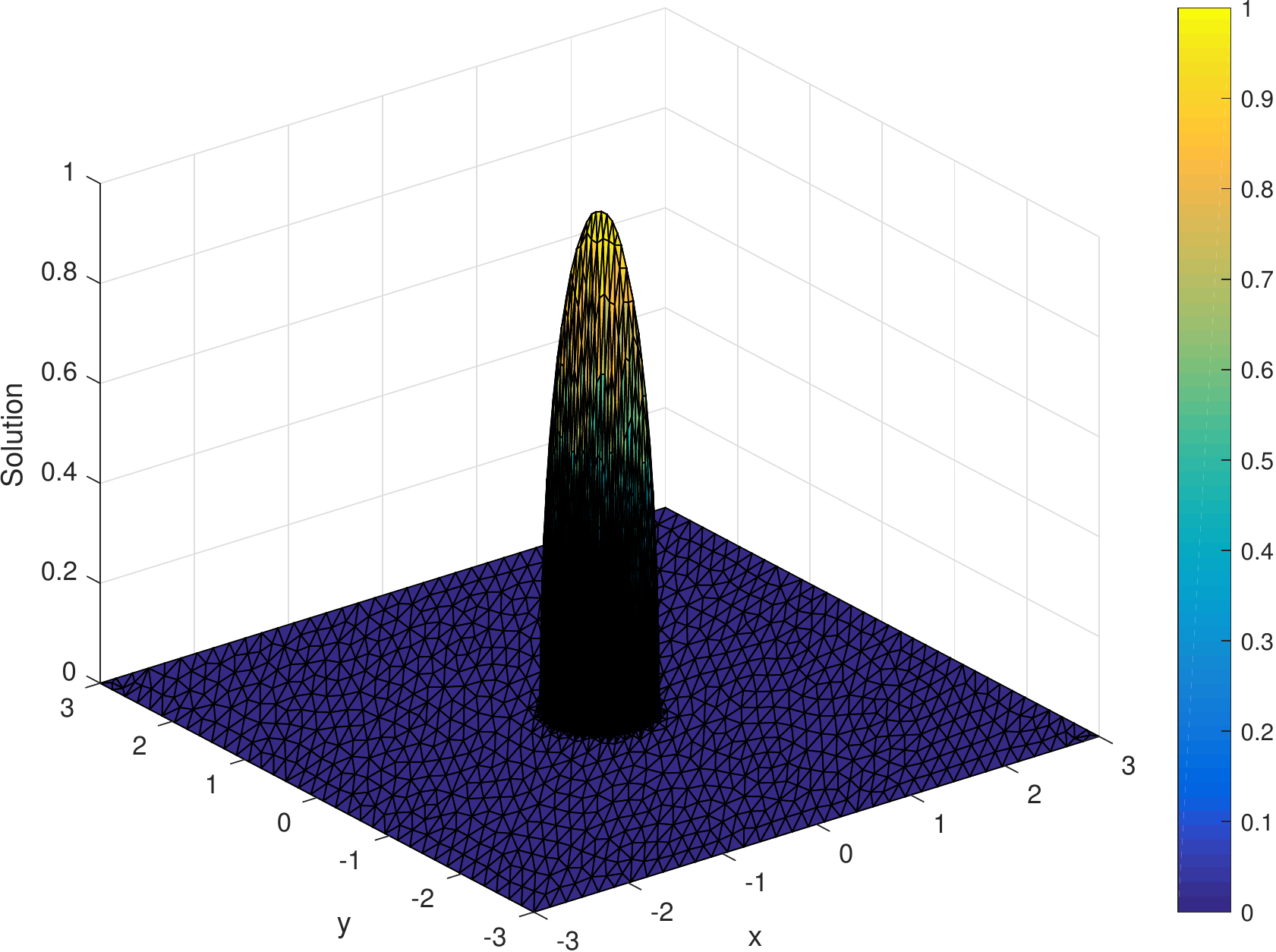}
\centerline{(a): $t=t_0$ }
\end{minipage}
\begin{minipage}{2.5in}
\includegraphics[width=2.5in]{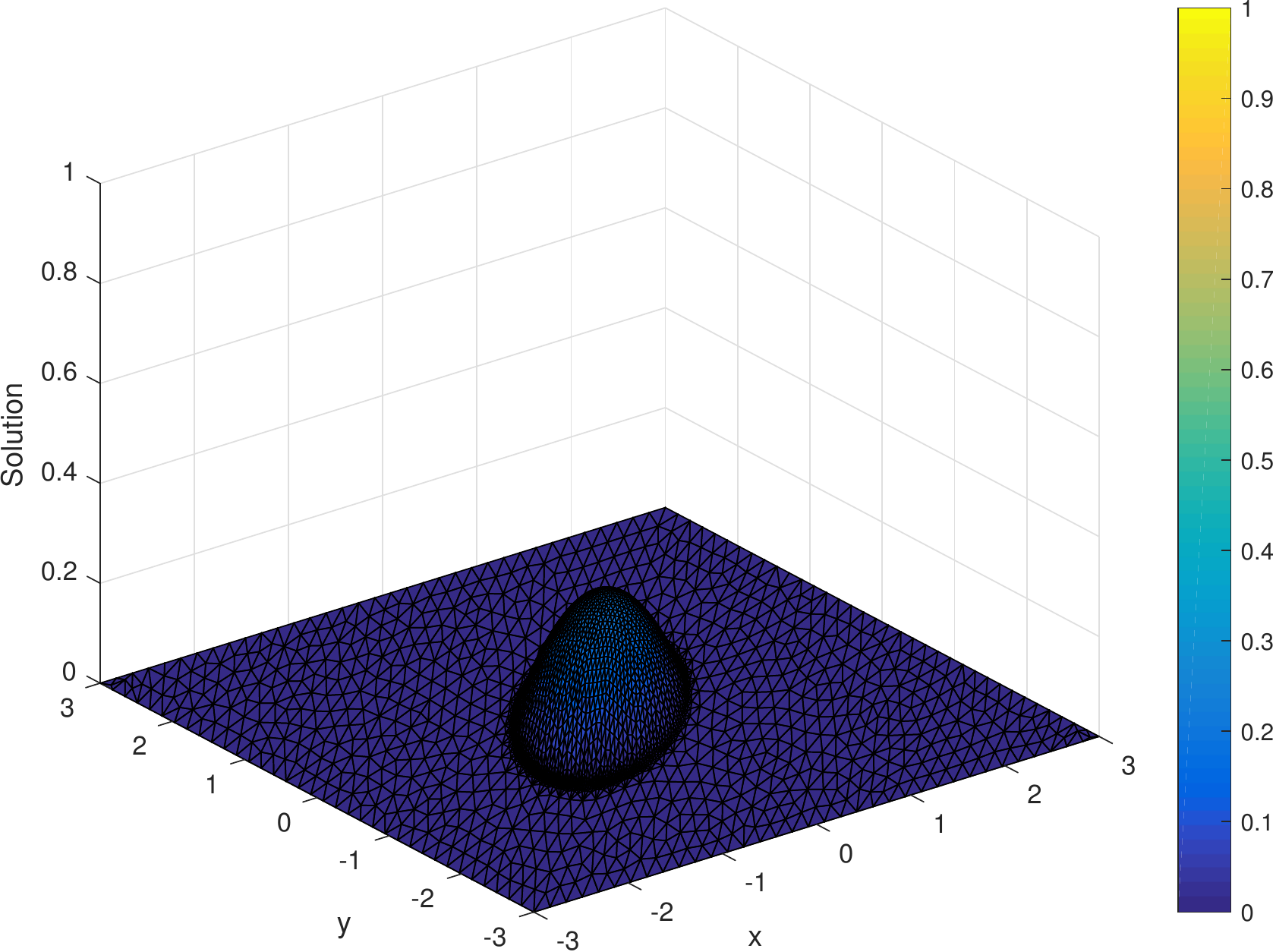}
\centerline{(b): $t=0.2$}
\end{minipage}
}
\caption{Example \ref{ex2}. Initial and final solutions, $m=1$.}
\label{ex2-soln}
\end{figure}

For investigation purpose, the results for $m=2$ using the same initial solutions are presented in Fig. \ref{ex2-meshb}. The observations are similar to the results for $m=1$ except that the free boundary has not moved as much as that for $m=1$ during the same time period. 
\begin{figure}[!thb]
\centering
\hbox{
\begin{minipage}{2.5in}
\includegraphics[width=2.5in]{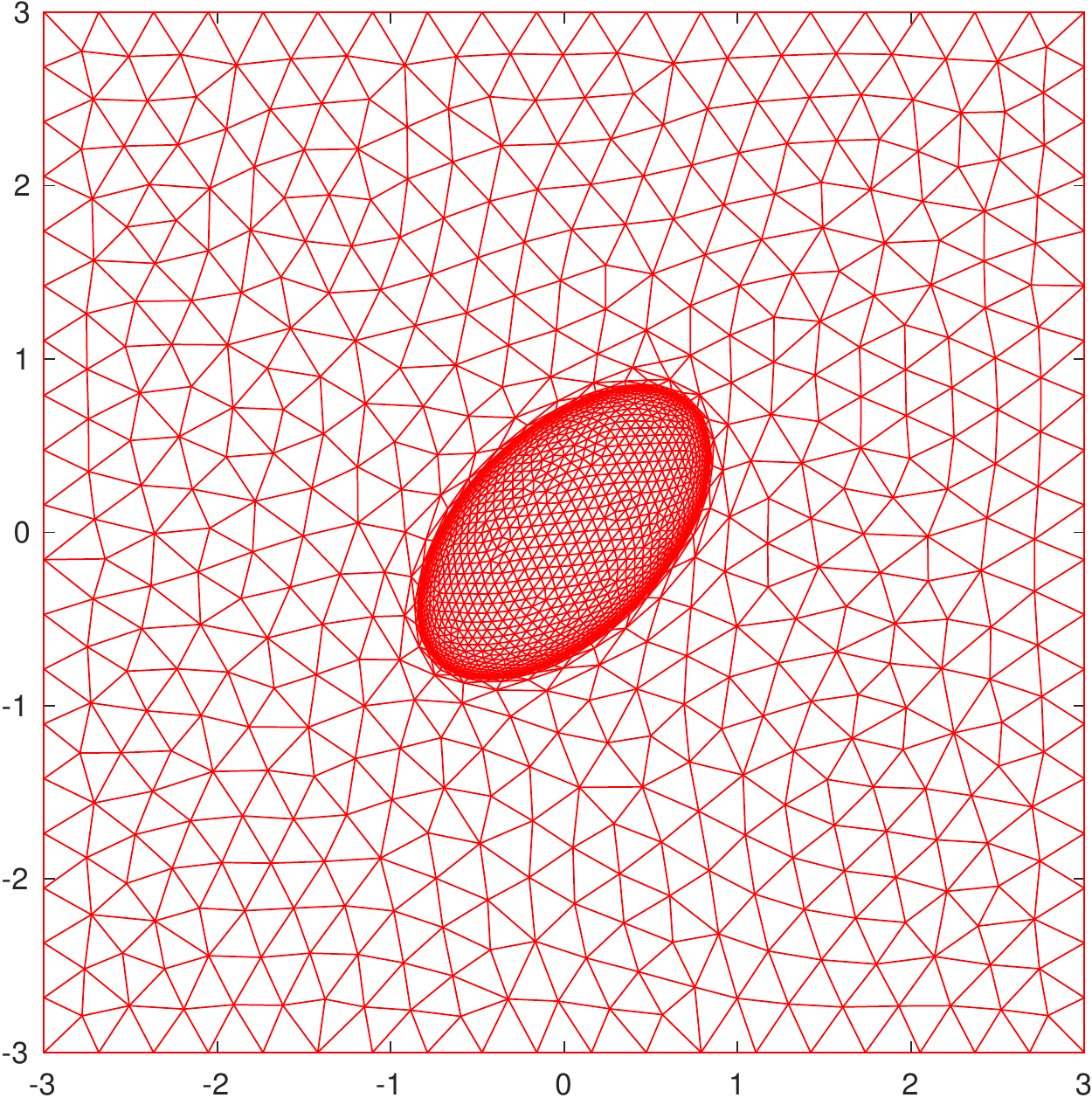}
\centerline{(a): $\M_{adap}$ mesh}
\end{minipage}
\begin{minipage}{2.5in}
\includegraphics[width=2.5in]{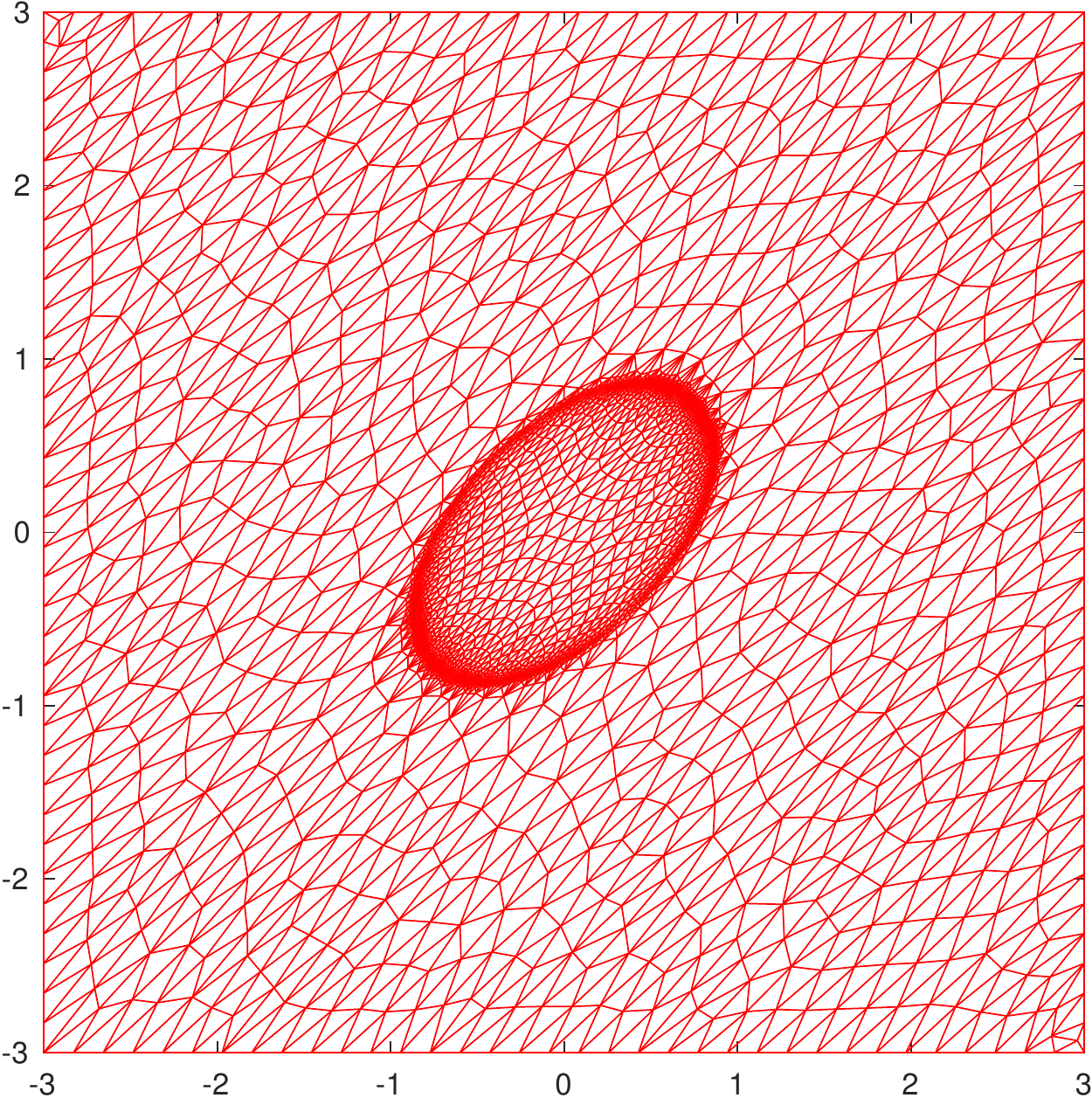}
\centerline{(b): $\M_{DMP+adap}$ mesh}
\end{minipage}
}
\caption{Example \ref{ex2}. $\M_{adap}$ and $\M_{DMP+adap}$ meshes at $t=0.2$, $m=2$.}
\label{ex2-meshb}
\end{figure}

\end{exam}

\vspace{10pt}

\begin{exam}
\label{ex3}
In the third example, we choose $m=6$ and consider an initial solution with two-isolated support in domain $\Omega=[-3,3]^2$. We also consider heterogeneous anisotropic diffusion. The initial solution is defined at $t_0=0$ as
\begin{equation}
	u_0 = \begin{cases}
		1, & \V{x} \in \Omega_1=(1,2) \times (0,1), \\
		1, & \V{x} \in \Omega_2=(-2,-1) \times (0.5,1.5) \\
		0, & \V{x} \in \Omega \backslash (\Omega_1 \cup \Omega_2).
	\end{cases}
\end{equation}
The end of time evolution is chosen as $T=120$. The diffusion matrix $\D=\D(\V{x})$ is taken as 
\begin{equation}
\D(\V{x}) = \left (\begin{array}{cc} \cos\theta & -\sin\theta \\ \sin\theta & \cos\theta \end{array} \right )
 \left (\begin{array}{cc} 50 & 0 \\ 0 & 1 \end{array} \right )
\left (\begin{array}{cc} \cos\theta & \sin\theta \\ -\sin\theta & \cos\theta \end{array} \right ),
\label{D-1}
\end{equation}
where $\theta = \theta(x,y) = \pi \sin(0.2x) \cos(0.1y)$ is the angle between the direction of the principle eigenvector of $\D$ and the positive $x$-axis. With this choice of $\theta$, the primary diffusion direction changes at different locations. 

The initial and final solutions are shown in Fig. \ref{ex3-soln}. The initial $\M_{adap}$ and $\M_{DMP+adap}$ meshes are displayed in Fig. \ref{ex3-mesh-ini}. As can be seen, the elements of the $\M_{DMP+adap}$ mesh are aligned along the principle diffusion directions at different places.

\begin{figure}[!thb]
\centering
\hbox{
\begin{minipage}{2.5in}
\includegraphics[width=2.5in]{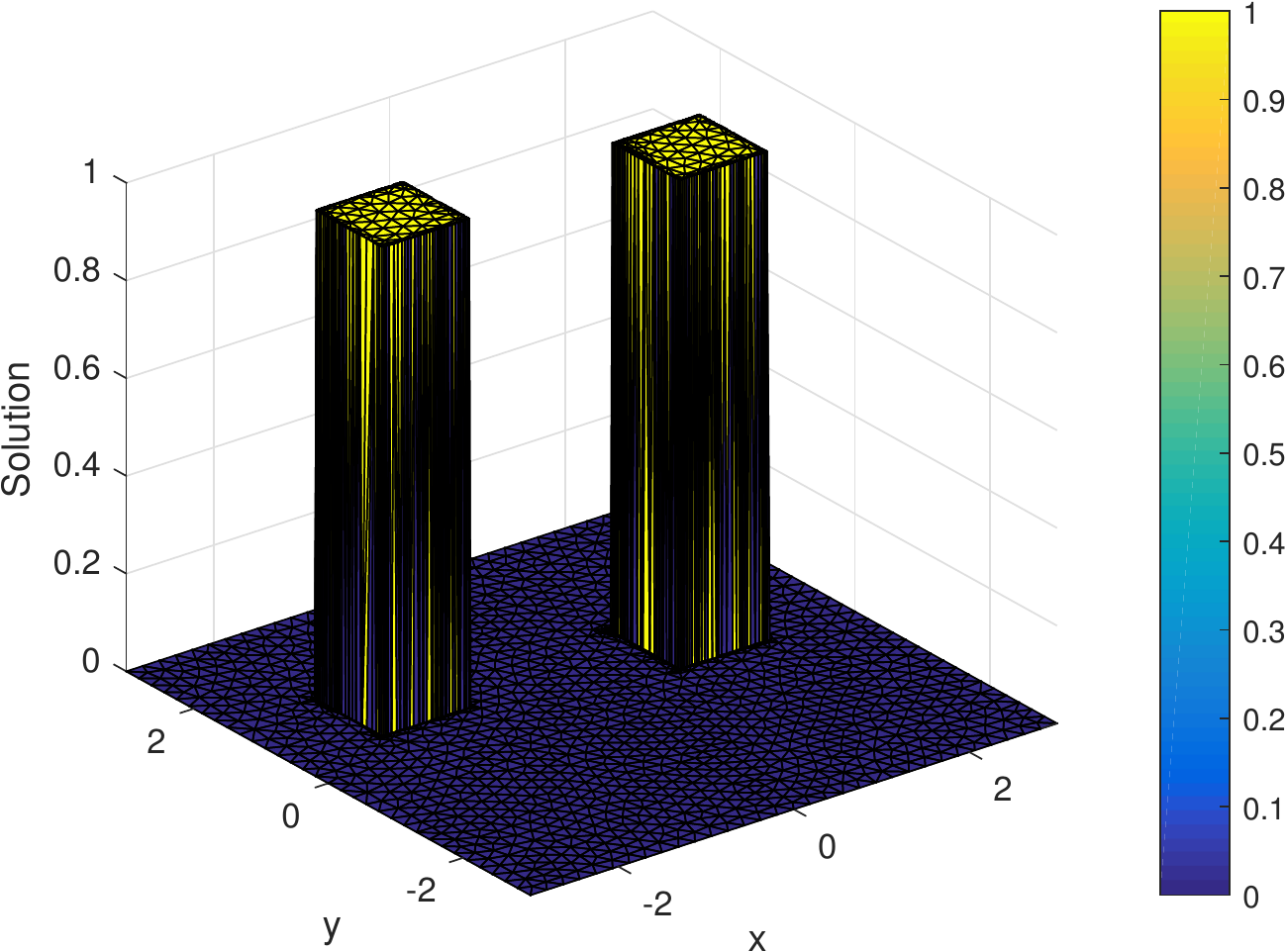}
\centerline{(a): $t=0$ }
\end{minipage}
\begin{minipage}{2.5in}
\includegraphics[width=2.5in]{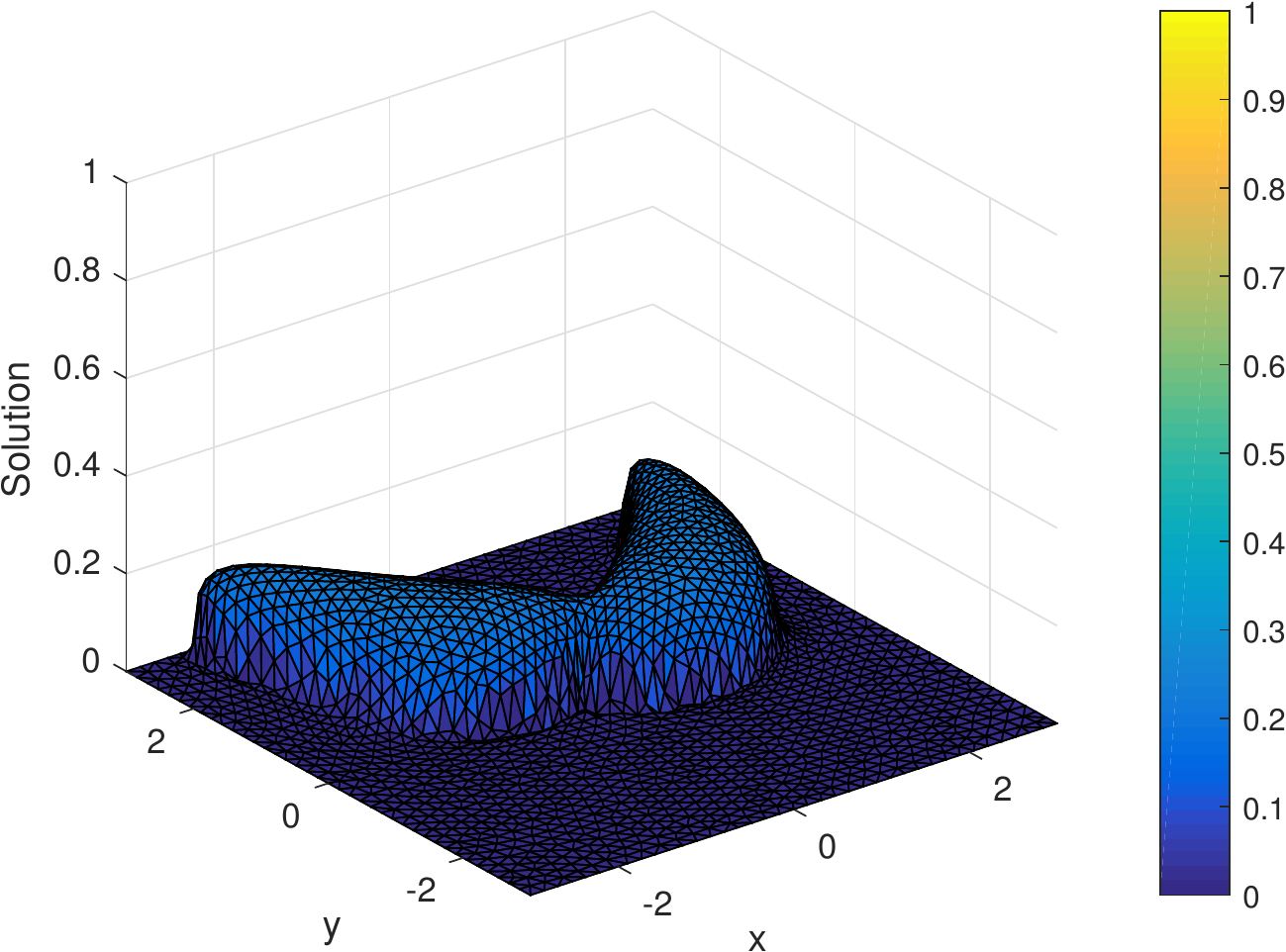}
\centerline{(b): $t=120$}
\end{minipage}
}
\caption{Example \ref{ex3}. Initial and final solutions, $m=6$.}
\label{ex3-soln}
\end{figure}

\begin{figure}[!thb]
\centering
\hbox{
\begin{minipage}{2.5in}
\includegraphics[width=2.5in]{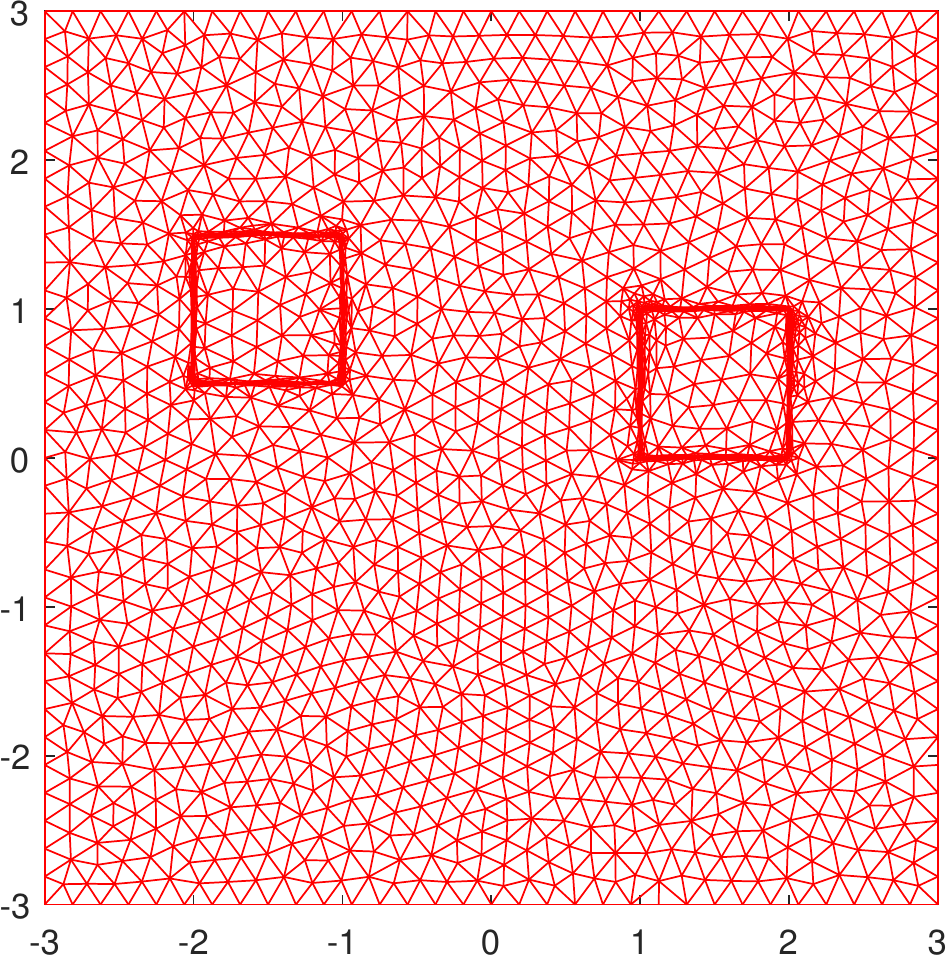}
\centerline{(a): $\M_{adap}$ mesh}
\end{minipage}
\begin{minipage}{2.5in}
\includegraphics[width=2.5in]{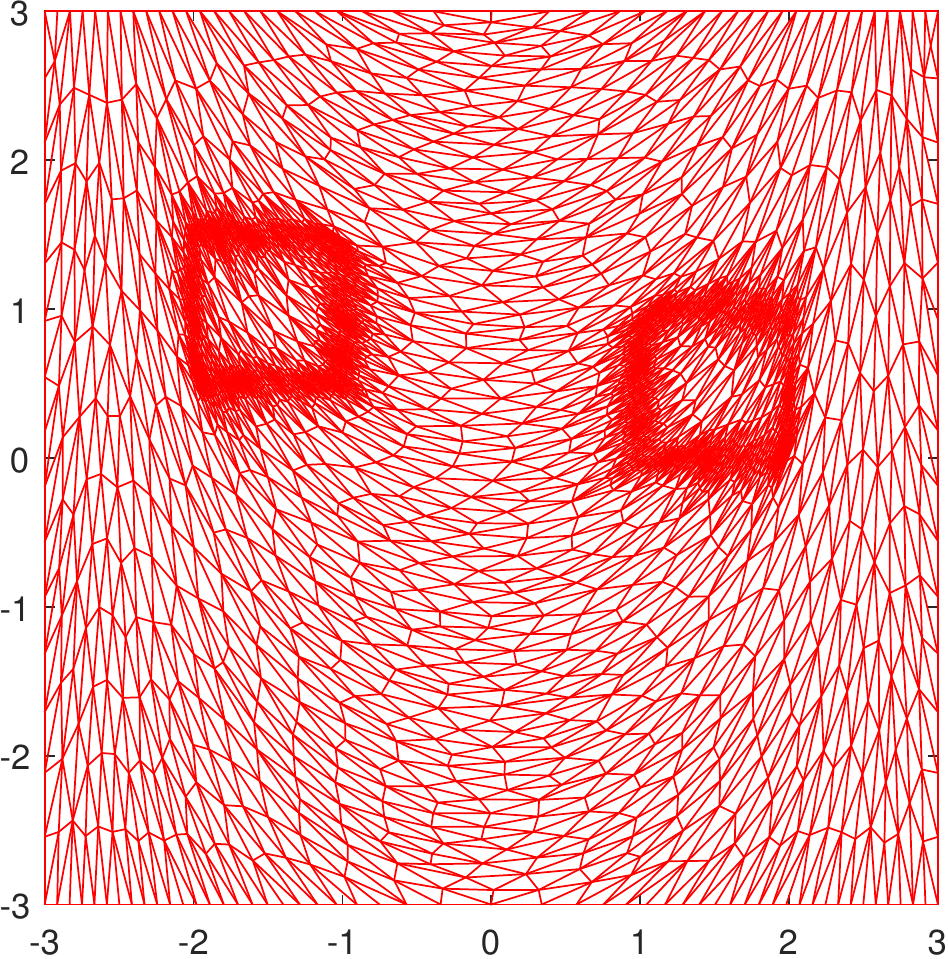}
\centerline{(b): $\M_{DMP+adap}$ mesh}
\end{minipage}
}
\caption{Example \ref{ex3}. Initial $\M_{adap}$ and $\M_{DMP+adap}$ meshes, $m=6$.}
\label{ex3-mesh-ini}
\end{figure}

Fig. \ref{ex3-meshb} and Fig. \ref{ex3-solnb} show the $\M_{adap}$ mesh and the corresponding numerical solutions at different times, respectively. It is observed that the two free boundaries start merging at around $t=60$. 

\begin{figure}[!thb]
\centering
\hbox{
\begin{minipage}{2.5in}
\includegraphics[width=2.5in]{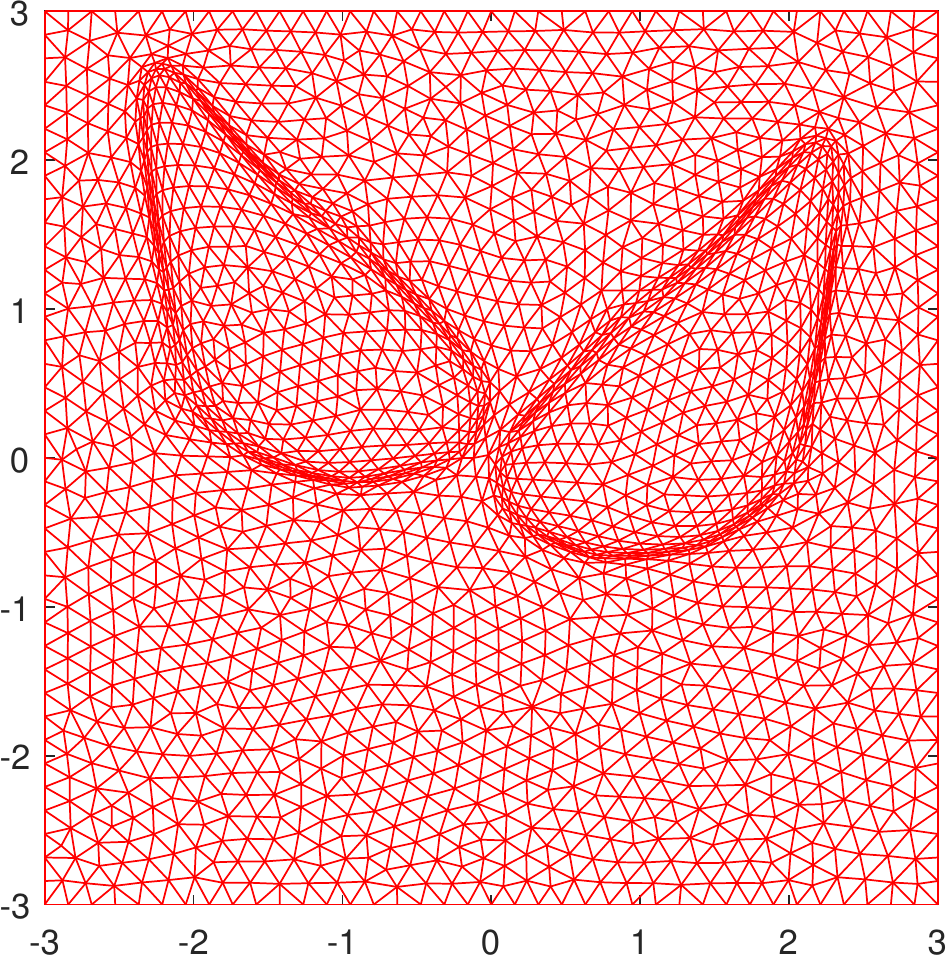}
\centerline{(a): $t=30$}
\end{minipage}
\begin{minipage}{2.5in}
\includegraphics[width=2.5in]{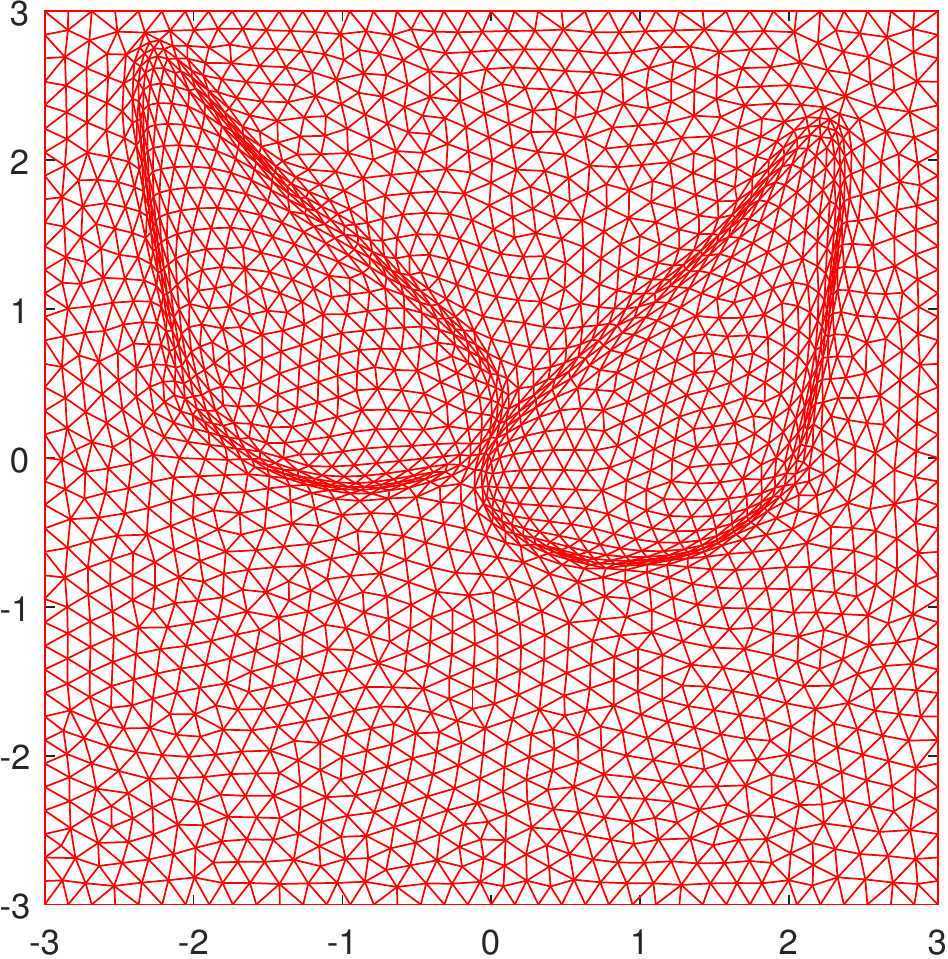}
\centerline{(b): $t=60$}
\end{minipage}
}
\vspace{5mm}
\hbox{
\begin{minipage}{2.5in}
\includegraphics[width=2.5in]{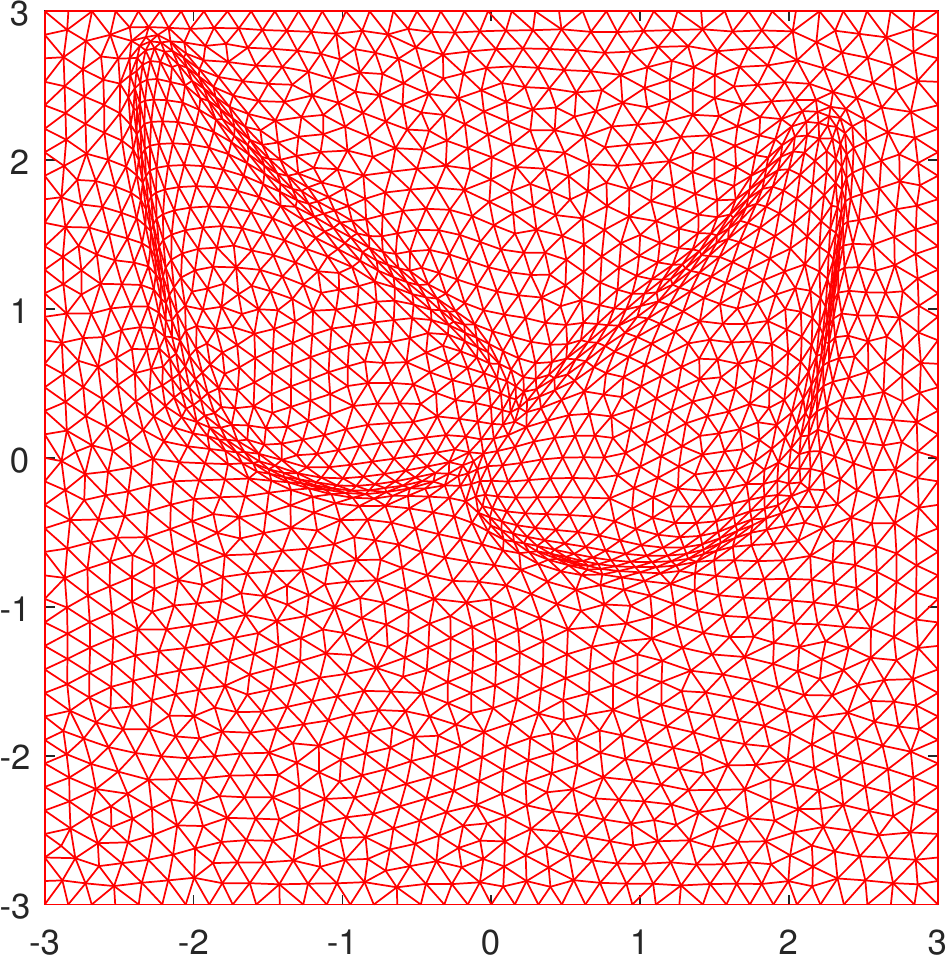}
\centerline{(c):  $t=90$}
\end{minipage}
\begin{minipage}{2.5in}
\includegraphics[width=2.5in]{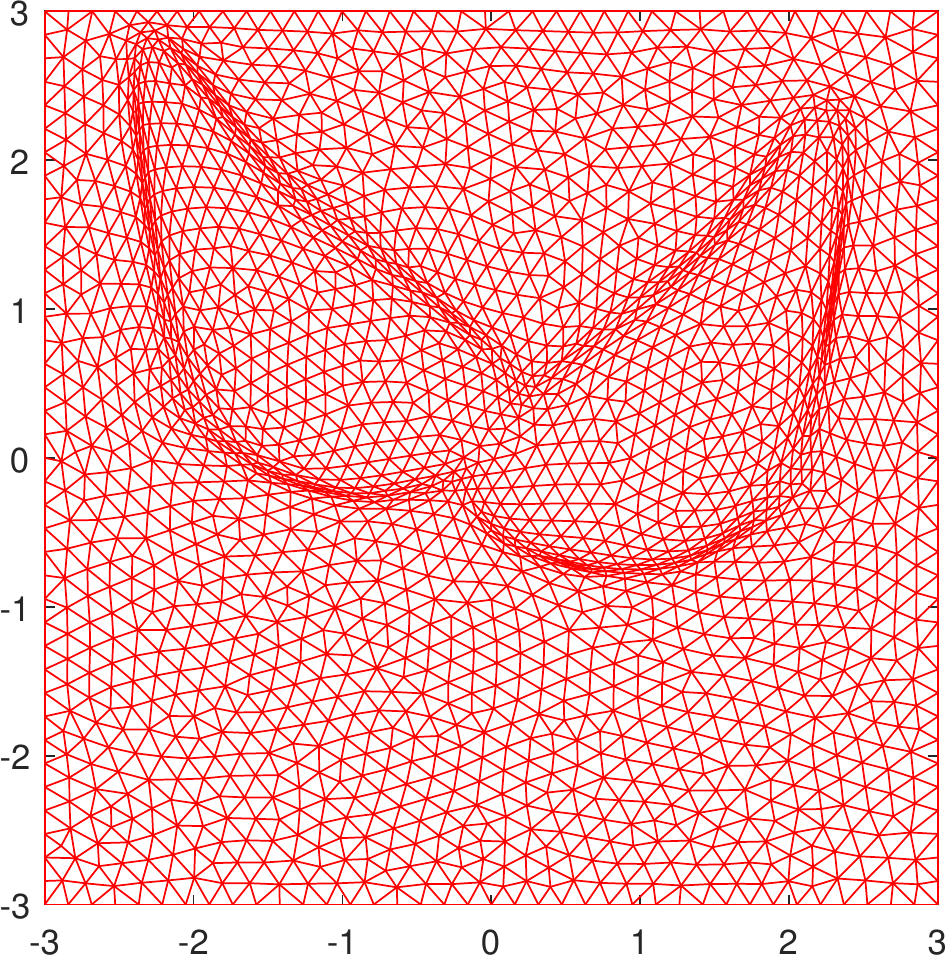}
\centerline{(d): $t=120$}
\end{minipage}
}
\caption{Example \ref{ex3}. $\M_{adap}$ meshes at different times, $m=6$.}
\label{ex3-meshb}
\end{figure}

\begin{figure}[!thb]
\centering
\hbox{
\begin{minipage}{2.5in}
\includegraphics[width=2.5in]{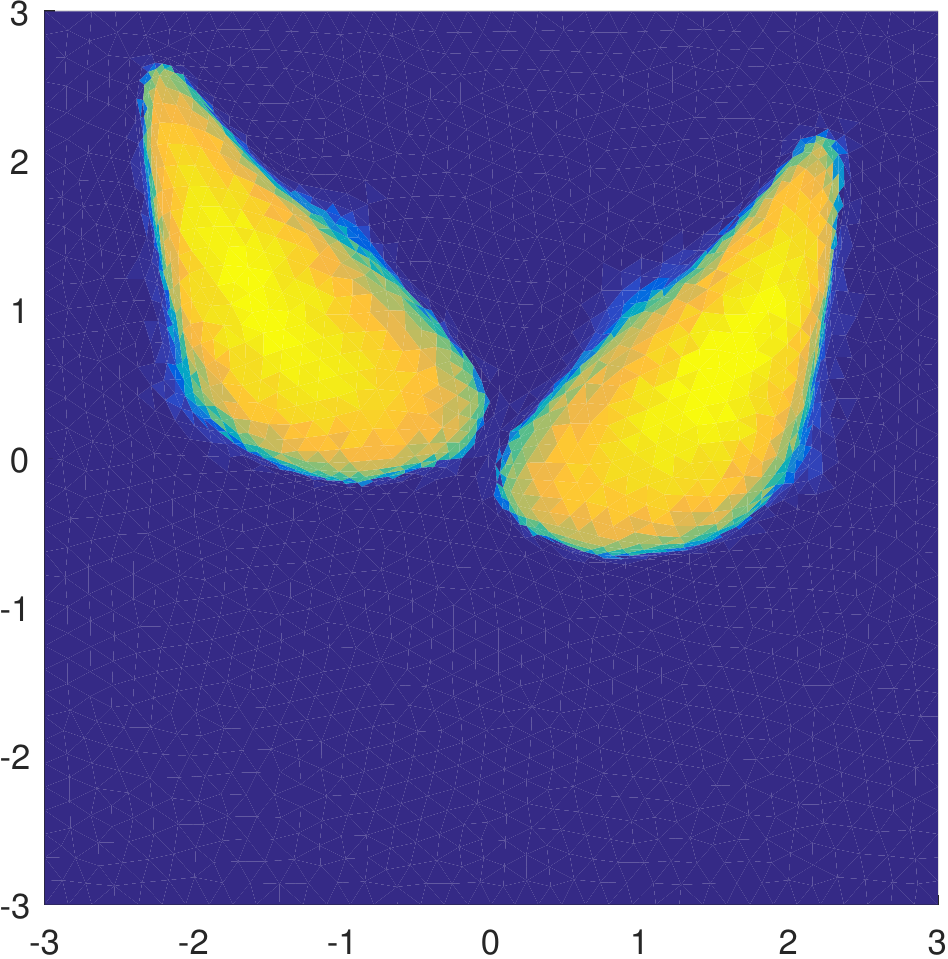}
\centerline{(a): $t=30$}
\end{minipage}
\begin{minipage}{2.5in}
\includegraphics[width=2.5in]{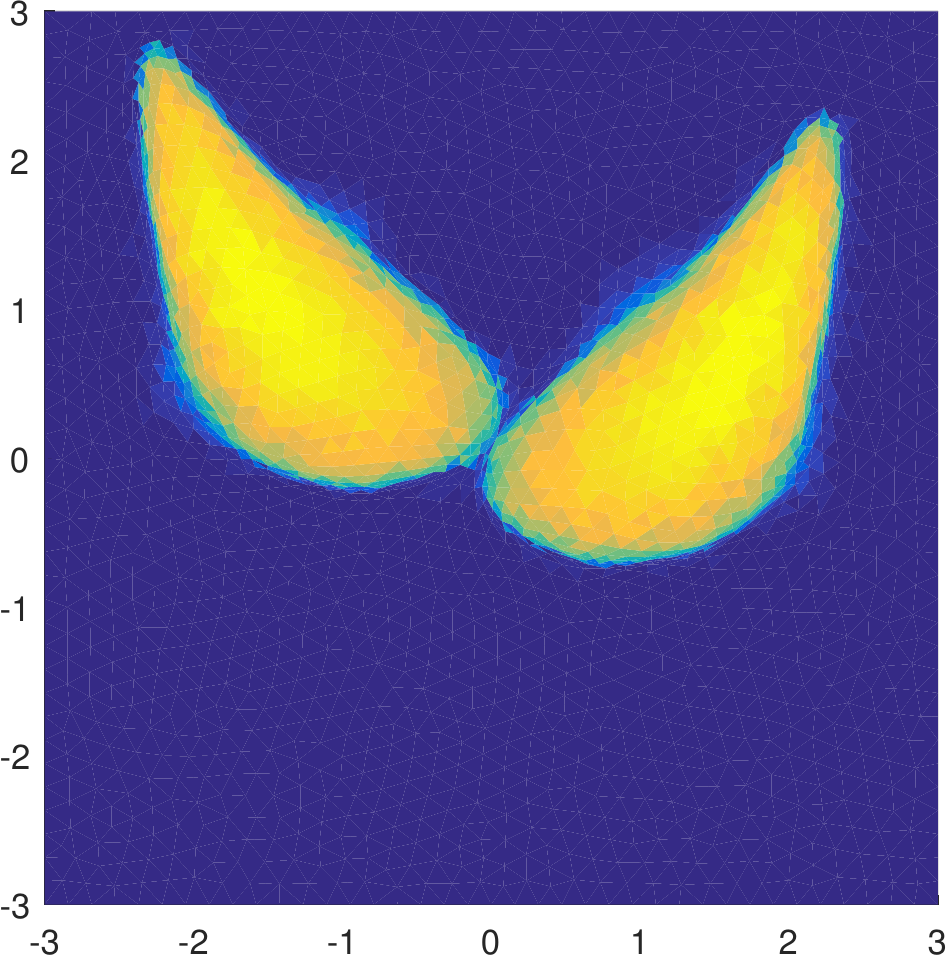}
\centerline{(b): $t=60$}
\end{minipage}
}
\vspace{5mm}
\hbox{
\begin{minipage}{2.5in}
\includegraphics[width=2.5in]{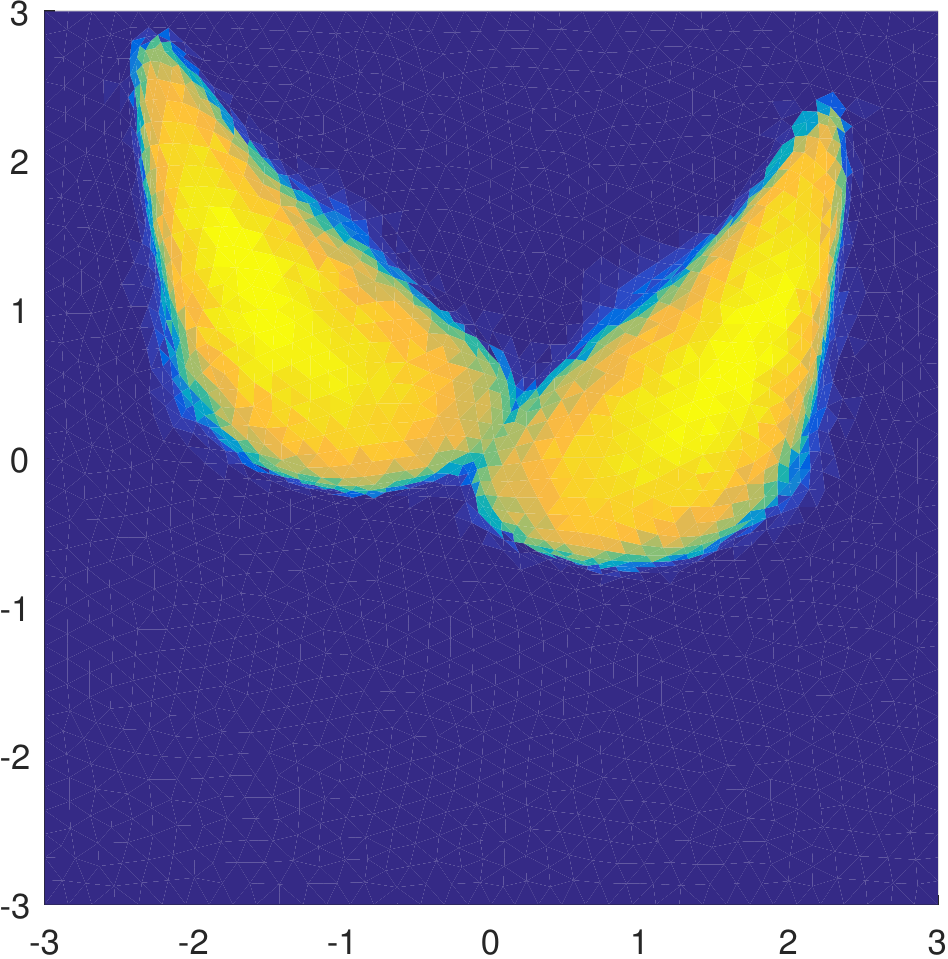}
\centerline{(c):  $t=90$}
\end{minipage}
\begin{minipage}{2.5in}
\includegraphics[width=2.5in]{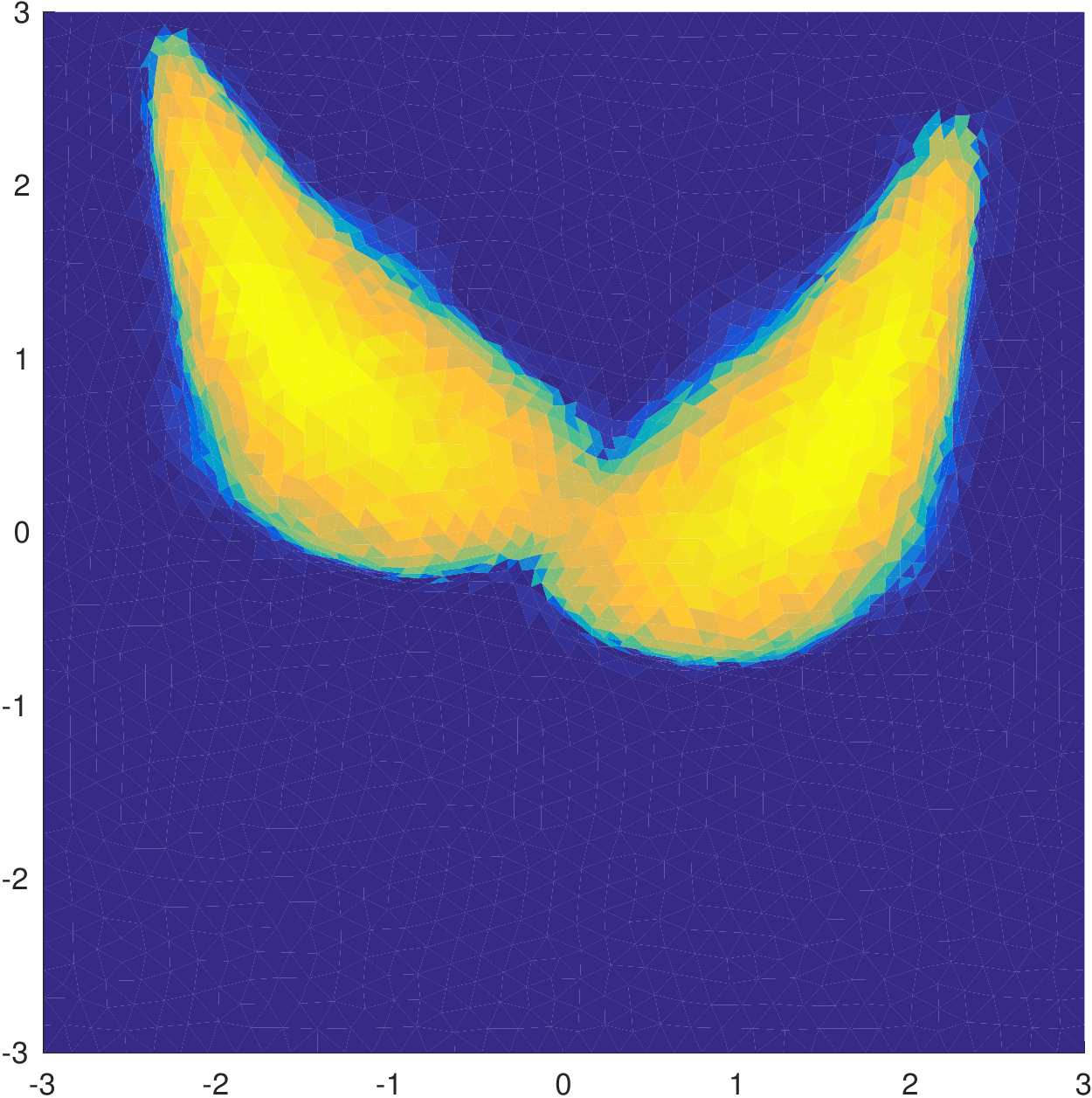}
\centerline{(d): $t=120$}
\end{minipage}
}
\caption{Example \ref{ex3}. Numerical solutions using $\M_{adap}$ meshes at different times, $m=6$.}
\label{ex3-solnb}
\end{figure}

\end{exam}

\begin{rem}
	Due to the challenges of the anisotropic and nonlinear diffusion in the APME, some negative values may occur in the numerical solutions during the computations, which violates the discrete maximum principle (DMP). The cut-off method \cite{LHV13} is applied to force all numerical solutions to be nonnegative. Satisfaction of DMP for APME and the effect of cut-off are under investigation. 
\end{rem}

\section{Summary}
\label{Sec-summary}

In the previous sections, we have generalized the Porous Medium Equation (PME) \eqref{pme-1} or \eqref{pme-2} to Anisotropic Porous Medium Equation (APME) \eqref{apme} that takes into account the anisotropy and heterogeneity of the physical properties of the porous media such as permeability. A special exact solution for APME is developed in \eqref{soln-apme} based on the Barenblatt-Pattle solution \eqref{soln-pme} for PME. 

Meanwhile, anisotropic mesh adaptation technique has been applied to obtain the finite element solutions of APME, which have helped improving both accuracy and efficiency of the computation. For general quasi-uniform meshes, the convergence rate for the solution error can be at mostly first order for $m=1$. However, with adaptive $\M_{adap}$ mesh or $\M_{DMP+adap}$ mesh, we have attained second order convergence rate, as demonstrated by  Example \ref{ex1}. Furthermore, for a same metric tensor, better adaptation and smaller errors can be achieved by adjusting the regularization parameter $\alpha_h$ in the computation. 

Our result is comparable with those obtained using moving mesh methods in \cite{NH17} for PME but can achieve the same accuracy with less number of mesh elements. Different than the moving mesh method that keeps the connectivity of the mesh elements, the adaptive meshes, including $\M_{adap}$ and $\M_{DMP+adap}$, not only change the connectivity but also can change the number of elements. Therefore, we can start with a coarse initial uniform mesh and adapt it to concentrate more elements around the free boundaries. The adapted mesh is used as a better initial mesh for later computations.      

Numerical results also show different behavior between APME and PME, as demonstrated by Examples \ref{ex1} and \ref{ex2}. For PME, the diffusion is the same in all directions, however, for APME, the diffusion is more significant along the principle diffusion direction. The merger between isolated free boundaries also occurs along the principle diffusion directions as demonstrated by Example \ref{ex3}.  The challenges for general anisotropic diffusion problems also apply to APME such as satisfaction of maximum principle \cite{KSS09, LH10, LH13}. More investigation on properties of APME is needed.

\vspace{20pt}

\textbf{Acknowledgment.} This work was partially supported by the UMRB grant KDM56 (Li) from the University of Missouri Research Board.   


\end{document}